\documentclass{amsart}
\usepackage[foot]{amsaddr}
\usepackage{graphicx}
\usepackage{color}
\usepackage{amssymb}
\usepackage[numbers,sort&compress]{natbib}
\usepackage{hyperref}
\usepackage{doi}
\usepackage{enumerate}

\newcommand{\dotprod}[2]{\langle #1, #2 \rangle}
\newcommand{\ddotprod}[2]{\langle\!\langle #1, #2 \rangle\!\rangle}
\newcommand{\argmin}{\operatornamewithlimits{argmin}}
\newcommand{\argmax}{\operatornamewithlimits{argmax}}

\def\RR{\mathbb{R}}  
\def\EE{\mathbb{E}}\def\PP{\mathbb{P}}

\def\Xint#1{\mathchoice
{\XXint\displaystyle\textstyle{#1}}%
{\XXint\textstyle\scriptstyle{#1}}%
{\XXint\scriptstyle\scriptscriptstyle{#1}}%
{\XXint\scriptscriptstyle%
\scriptscriptstyle{#1}}%
\!\int}
\def\XXint#1#2#3{{\setbox0=\hbox{$#1{#2#3}{%
\int}$ }
\vcenter{\hbox{$#2#3$ }}\kern-.6\wd0}}

\def\dashint{\Xint-}

\renewcommand{\pmb}{} 

\begin{document}

\title[Long Term Effects of Small Random Perturbations]{Long Term
  Effects of Small Random Perturbations on Dynamical Systems:
  Theoretical and Computational Tools}

\author{Tobias Grafke$^1$, Tobias Sch\"afer$^2$, and Eric
  Vanden-Eijnden$^1$}

\email{grafke@cims.nyu.edu, tobias.schaefer@csi.cuny.edu, eve2@cims.nyu.edu}
\address{$^1$Courant Institute, 
  New York University, 
  251 Mercer Street, New York, 
  NY 10012, USA}
\address{$^2$Department of Mathematics, 
  College of Staten Island 1S-215, 2800 Victory Blvd., 
  Staten Island, New York 10314}
\date{\today}

\begin{abstract}
  Small random perturbations may have a dramatic impact on the
  long time evolution of dynamical systems, and large deviation theory
  is often the right theoretical framework to understand these
  effects. At the core of the theory lies the minimization of an
  action functional, which in many cases of interest has to be
  computed by numerical means. Here we review the theoretical and
  computational aspects behind these calculations, and propose an
  algorithm that simplifies the geometric minimum action method to
  minimize the action in the space of arc-length parametrized
  curves. We then illustrate this algorithm's capabilities by applying
  it to various examples from material sciences, fluid dynamics,
  atmosphere/ocean sciences, and reaction kinetics. In terms of
  models, these examples involve stochastic (ordinary or partial)
  differential equations with multiplicative or degenerate noise,
  Markov jump processes, and systems with fast and slow degrees of
  freedom, which all violate detailed balance, so that simpler
  computational methods are not applicable.
\end{abstract}

\maketitle

\section{Introduction}
\label{sec:intro}

Small random perturbations often have a lasting effect on the
long-time evolution of dynamical systems. For example, they give rise
to transitions between otherwise stable equilibria, a phenomenon
referred to as metastability that is observed in a wide variety of
contexts, e.g.  phase separation, population dynamics, chemical
reactions, climate regimes, neuroscience, or fluid dynamics. Since the
time-scale over which these transition events occurs is typically
exponentially large in some control parameter (for example the noise
amplitude), a brute-force simulation approach to compute these events
quickly becomes infeasible. Fortunately, it is possible to exploit the
fact that the mechanism of these transitions is often predictable when
the random perturbations have small amplitude: with high probability
the transitions occur by their path of maximum likelihood\index{path
  of maximum likelihood} (PML), and knowledge of this PML also permits
to estimate their rate. This is the essence of large deviation theory
(LDT) \cite{freidlin-wentzell:2012}, which applies in a wide variety
of contexts. For example, systems whose evolution is governed by a
stochastic (ordinary or partial) differential equation driven by a
small noise or by a Markov jump process in which jumps occur often but
lead to small changes of the system state, or slow/fast systems in
which the fast variables are randomly driven and the slow ones feel
these perturbations through the effect fast variables only, all fit
within the framework of LDT. Note that, typically, the dynamics of
these systems fail to exhibit microscopic reversibility (detailed
balance)\index{detailed balance} and the transitions therefore occur
out-of-equilibrium. Nevertheless, LDT still applies.

LDT also indicates that the PML is computable as the minimizer of a
specific objective function (action): the large deviation rate
function of the problem at hand. 
This is a non-trivial numerical optimization problem which calls for
tailor-made techniques for its solution.  Here we will focus on one
such technique, the geometric minimum action method (gMAM,
\cite{vanden-eijnden-heymann:2008, heymann-vanden-eijnden:2008-b,
  vanden-eijnden-heymann:2008}), which builds the minimum action
method and its variants \cite{e-ren-vanden-eijnden:2004,%
  zhou-ren-e:2008, wan:2011}, and was designed to perform the action
minimization over both the transition path location and its
duration. This computation gives the so-called quasipotential, whose
role is key to understand the long time effect of the random
perturbations on the system, including the mechanism of transitions
events induced by these perturbations. Our purpose here is twofold.
First, we would like to briefly review the theoretical aspects behind
LDT that led to the rate function minimization problem and, in
particular, to the geometric variant of it that is central in
gMAM. Second, we would like to discuss in some details the
computational issues this minimization entails, and remedy a drawback
of gMAM, namely its somewhat complicated descent step that requires
higher order derivatives of the large deviation Hamiltonian. Here, we
propose a simpler algorithm, minimizing the geometric action
functional, but requiring only first order derivatives of the
Hamiltonian. The power of this algorithm is then illustrated via
applications to a selection of problems:
\begin{enumerate}
\item the Maier-Stein model, which is a toy non-gradient stochastic
  ordinary differential equation that breaks detailed balance;
\item a stochastic Allen-Cahn/Cahn-Hilliard partial differential
  equation motivated by population dynamics;
\item the stochastic Burgers-Huxley PDE, related to fluid dynamics and
  neuroscience;
\item Egger's and Charney-DeVore equations, introduced as climate
  models displaying noise-induced transitions between metastable
  regimes;
\item a generalized voter/Ising model with multiplicative noise;
\item metastable networks of chemical reaction equations and
  reaction-diffusion equations;
\item a fast/slow system displaying transitions of the slow variables
  induced by the effects of the fast ones.
\end{enumerate}

The remainder of this paper is organized as follows.  In
Sect.~\ref{ssec:ldt} we briefly review the key concepts of LDT that we
will use (Sect.~\ref{sec:key}) and give a geometrical point of view of
the theory that led to the action used in gMAM (sec.~\ref{ssec:gmam}).
In Sect.~\ref{sec:algorithm} we discuss the numerical aspects related
to the minimization of the geometric action, propose a simplified
algorithm to perform this calculation, and compare it to existing
algorithms. We also discuss further simplifications of the algorithm
that apply in regularly occurring special cases, such as additive or
multiplicative Gaussian noise. Finally, in Sect.~\ref{sec:applications}
we present the applications listed above.

\section{Freidlin-Wentzell large deviation theory (LDT)}
\label{ssec:ldt}

Here we first give a brief overview of
LDT~\cite{freidlin-wentzell:2012}, focusing mainly on stochastic
differential equations (SDEs) for simplicity, but indicating also how
the theory can be extended to other models, such as Markov jump
processes or fast/slow systems. Then we discuss the geometric
reformulation of the action minimization problem that is used in gMAM.

\subsection{Some key concepts in LDT}
\label{sec:key}

Consider the following SDE for $X\in \RR^n$
\begin{equation}
  \label{eq:sde}
  dX = b(X) dt + \sqrt{\epsilon} \sigma(X) dW\,,
\end{equation}
where $b: \RR^n \to\RR^n$ denotes the drift term, $W$ is a standard
Wiener process in $\RR^n$, $\sigma: \RR^n \to \RR^n \times \RR^n$ is
related to the diffusion tensor via $a(x) = (\sigma
\sigma^\dagger)(x)$, and $\epsilon>0$ is a parameter measuring the
noise amplitude.  Suppose that we want to estimate the probability of
an event, such as finding the solution in a set $B\subset \RR^n$ at
time $T$ given that it started at $X(0)=x$ at time $t=0$. LDT
indicates that, in the limit as $\epsilon\to0$, this probability can
be estimated via a minimization problem\index{large deviation theory}:
\begin{equation}
  \label{eq:1}
  \PP^x \left(X(T)\in B\right) \asymp \exp\left(-\epsilon^{-1} \min _{\phi\in
      \mathcal {C}} S_T(\phi)  \right)\,.
\end{equation}
Here $\asymp$ denotes log-asymptotic equivalence\index{log-asymptotic
  equivalence} (i.e. the ratio of the logarithms of both sides tends
to 1 as $\epsilon\to0$), the minimum is taken over the set
$\mathcal{C} = \{ \phi\in C([0,T],\RR^n): \phi(0)=x,\phi(T)\in B\}$,
and we defined the action functional\index{action functional}
\begin{equation}
  \label{eq:action}
  S_T(\phi) = 
  \begin{cases}
    \int_0^T L(\phi, \dot \phi)\,dt & \text{if the integral converges}\\
    \infty & \text{otherwise.}
  \end{cases}
\end{equation}
Here 
\begin{equation}
  \label{eq:lagrange}
  L(\phi,\dot \phi) = \tfrac12\dotprod{\dot \phi - b(\phi)}
  {\left(a(\phi)\right)^{-1}(\dot \phi - b(\phi))}\,,
\end{equation}
where we assumed for simplicity that $a(\phi)$ is invertible -- (this
assumption will be relaxed below) and $\dotprod{\cdot}{\cdot}$ denotes
the Euclidean inner product in $\RR^n$. LDT also indicates that, as
$\epsilon \to0$, when the event occurs, it does so with $X$ being
arbitrarily close to the minimizer
\begin{equation}
  \label{eq:2}
  \phi_* = \mathop{\text{argmin }} _{\phi\in
      \mathcal {C}} S_T(\phi)
\end{equation}
in the sense that
\begin{equation*}
  \label{eq:8}
  \forall \delta>0: \qquad \lim_{\epsilon\to0}\PP^x \Bigl(\,\sup_{0\le t\le T}|X(t) -
  \phi_*(t) \,|\, < \delta \Big| X(T) \in B\Bigr) = 1
\end{equation*}
Thus, from a computational viewpoint, the main question becomes how to
perform the minimization in~\eqref{eq:2}. Note that, if we define the
Hamiltonian\index{Hamiltonian} associated with the
Lagrangian\index{Lagrangian}~\eqref{eq:lagrange}
\begin{equation}
  \label{eq:H}
  H(\phi,\theta) = \dotprod{b(\phi)}{\theta}+\tfrac12\dotprod{\theta}{a(\phi) \theta}
\end{equation}
such that
\begin{equation}
  \label{eq:legendre}
  L(\phi,\dot\phi) = \sup_{\theta} \left(\dotprod{\dot \phi}{\theta} 
    - H(\phi,\theta)\right),
\end{equation}
this minimization reduces to the solution of Hamilton's equations of
motion\index{equations of motion},
\begin{equation}
  \label{eq:eqnmotion}
  \begin{cases}
    \dot \phi = H_\theta(\phi,\theta) = b(\phi) + a(\phi) \theta\\
    \dot \theta = -H_\phi(\phi,\theta) = -(b_\phi(\phi))^T \theta 
    + \tfrac12 \dotprod{\theta}{a_\phi(\phi) \theta}\,,
  \end{cases}
\end{equation}
where subscripts denote differentiation and we use the convention
$(b_\phi)_{ij} = \partial b_i/\partial \phi_j$. What makes the problem
nonstandard, however, is the fact that these equations must be solved
as a boundary value problem, with $\phi(0)=x$ and $\phi(T) = y\in
B$. We will come back to this issue below.

If the minimum of the action in~\eqref{eq:1} is nonzero, this equation
indicates that the probability of finding the solution in~$B$ at
time~$T$ is exponentially small in $\epsilon$, i.e. it is a rare
event. This is typically the case if one considers events that occur
on a finite time interval, $T<\infty$ fixed. LDT, however, also
permits to analyze the effects of the perturbations over an infinite
time span, in which case they become ubiquitous. In this context, the
central object in LDT is the quasipotential\index{quasipotential}
defined as
\begin{equation}
  \label{eq:minimize-T}
  V(x,y) = \inf_{T>0} \min_{\phi\in \mathcal{C}_{x,y}} S_T(\phi)\,,
\end{equation}
where $\mathcal{C}_{x,y} = \{ \phi\in
C([0,T],\RR^n):\phi(0)=x,\phi(T)=y\}$. The quasipotential permits to
answer several questions about the long time behavior of the
system. For example, if we assume that the deterministic equation
associated with~\eqref{eq:sde}, $\dot X= b(X)$, possesses a single
stable fixed point, $x_a$, and that \eqref{eq:sde} admits a unique
invariant distribution\index{invariant distribution}, the density
associated with this distribution can estimated as $\epsilon\to0$ as
\begin{equation}
  \label{eq:5}
  \rho(x) \asymp \exp\left(- \epsilon^{-1} V(x_a,x)\right)\,.
\end{equation}
Similarly, if $\dot X = b(X)$ possesses two stable fixed points, $x_a$
and $x_b$, whose basins of attraction have a common boundary, we can
estimate the mean first passage time\index{mean first passage time}
the system takes to travel for one fixed point to the other as
\begin{equation}
  \label{eq:4}
  \EE \tau_{a\to b} \asymp \exp\left(\epsilon^{-1} V(x_a,x_b)\right)\,,
\end{equation}
where
\begin{equation}
  \label{eq:3}
  \tau_{a\to b} = \inf \{ t: X(t)\in B_\delta (x_b), X(0) = x_a\}\,,
\end{equation}
in which $B_\delta(x_b)$ denotes the ball of radius $\delta$ around
$x_b$, with $\delta$ small enough so that this ball is contained in
the basin of attraction of $x_b$. In this set up, we can also estimate
the ratio of the stationary probabilities to find the system in the
basins of attraction of $x_a$ or $x_b$. Denoting these probabilities
by $p_a$ and $p_b$, respectively, we have
\begin{equation}
  \label{eq:6}
  \frac{p_a}{p_b} \asymp \frac{\EE \tau_{a\to b}}{\EE \tau_{b\to a}}
  \asymp
  \exp\left(\epsilon^{-1} (V(x_a,x_b)-V(x_b,x_a))\right)\,.
\end{equation}

These statements can be generalized to many other situations, e.g. if
$\dot X = b(X)$ possesses more than two stable fixed points, or
attracting structures that are more complicated than points, such as
limit cycles. They can also be generalized to dynamical systems other
than~\eqref{eq:sde}, e.g. if this equation is replaced by a stochastic
partial differential equation (SPDE), or for Markov jump processes in which
the jump rates are fast but lead to small changes of the system's
state~\cite{freidlin-wentzell:2012,shwartz-weiss:1995}, or in
slow/fast systems where the slow variables feels random perturbations
through the effect the fast variables have on
them~\cite{freidlin:1978, kifer:1992, veretennikov:2000, kifer:2004,
  bouchet-grafke-tangarife-etal:2016}. In all cases, LDT provides us
with an action functional like~\eqref{eq:action}, but in which the
Lagrangian is different from~\eqref{eq:lagrange} if the system's
dynamics is not governed by an S(P)DE. Typically, the theory yields an
expression for the Hamiltonian~\eqref{eq:H}, which may be
non-quadratic in the momenta, such that the Legendre transform
in~\eqref{eq:legendre} is not available analytically. This
\textit{per~se} is not an issue, since we can in principle minimize
the action by solving Hamilton's
equations~\eqref{eq:eqnmotion}. However, these calculations face two
difficulties. The first, already mentioned above, is
that~\eqref{eq:eqnmotion} must be solved as a boundary value
problem. The second, which is specific to the calculation of the
quasipotential in~\eqref{eq:minimize-T}, is that the time span over
which \eqref{eq:eqnmotion} are solved must be varied as well
since~\eqref{eq:minimize-T} involves a minimization over $T$, and
typically the minimum is reached as $T\to\infty$ (i.e. there is a
minimizing sequence but no minimizer) which complicates matters even
more. These issues motivate a geometric reformulation of the problem,
which was first proposed in \cite{heymann-vanden-eijnden:2008} and we
recall next.

\subsection{Geometric action functional}
\label{ssec:gmam}

As detailed in \cite{heymann-vanden-eijnden:2008} (see Proposition 2.1
in that paper), the quasipotential defined in~\eqref{eq:minimize-T}
can also be expressed as
\begin{equation}
  \label{eq:QP2}
  V(x,y) = \min_{\varphi\in \hat{\mathcal{C}}_{x,y}} \hat S(\varphi)\,,
\end{equation}
where $\hat{\mathcal{C}}_{x,y} = \{ \varphi\in
C([0,1],\RR^n):\varphi(0)=x,\varphi(1)=y\}$ and $\hat S(\varphi)$ is
the geometric action\index{geometric action functional} that can be defined in
the following equivalent ways:
\begin{subequations}
  \label{eq:gmam}
  \begin{align}
    \hat S(\varphi) &= \sup_{\vartheta: H(\varphi,\vartheta)=0} \int_0^1
                   \dotprod{\varphi'}{\vartheta} ds
                   \label{eq:gmam-2}\\
    \hat S(\varphi) &= \int_0^1 \dotprod{\varphi'}{\vartheta_*(\varphi,\varphi')} ds
                   \label{eq:gmam-3}\\
    \hat S(\varphi) &= \int_0^1 \frac1{\lambda(\varphi,\varphi')} L(\varphi,
                   \lambda \varphi') ds\,,
                   \label{eq:gmam-4}
  \end{align}
\end{subequations}
where $\vartheta_*(\varphi,\varphi')$ and $\lambda(\varphi,\varphi')$ are the
solutions to
\begin{equation}
  \label{eq:find-theta-gmam}
  H(\varphi,\vartheta_*(\varphi,\varphi')) = 0,\qquad 
  H_\vartheta(\varphi,\vartheta_*(\varphi,\varphi')) = \lambda(\varphi,\varphi') \varphi' 
  \quad\text{with } \lambda \ge 0\,.
\end{equation}
The action $\hat S(\varphi)$ has the property that its value is left
invariant by reparametrization of the path~$\varphi$, i.e. it is an
action on the space of continuous curves.  In particular, one is free
to choose arclength-parametrization for $\varphi$, e.g. $|\varphi'| =
1/L$ for $\int |\varphi'| \,ds = L$.  This also means that the
minimizer of~\eqref{eq:QP2} exists in more general cases (namely as
long as the path has finite length), which makes the minimization
problem easier to handle numerically, as shown next.

\section{Numerical minimization of the geometric action}
\label{sec:algorithm}

From~\eqref{eq:QP2}, we see that the calculation of the quasipotential
reduces to a minimization problem, whose Euler-Lagrange equation is
simply
\begin{equation}
  \label{eq:9}
  D_\phi\hat S(\varphi)=0\,,
\end{equation}
where $D_\varphi$ denotes the functional gradient with respect to
$\varphi$. The main issue then becomes how to find the solution
$\varphi_*$ to~\eqref{eq:9} that minimize the action
$\hat S(\varphi)$. In this section, we first briefly review how the
gMAM achieves this task. We will then introduce a simplified variant
of the gMAM algorithm that in its simplest form relies solely on first
order derivatives of the Hamiltonian. Subsequently, we also analyze
several special cases where the numerical treatment can be simplified
even further.

\subsection{Geometric minimum action method}

The starting point of gMAM\index{geometric minimum action method} is
the following expression involving $D_\phi \hat S(\varphi)$ that can
be calculated directly from formula~\eqref{eq:gmam-3} for the action
functional:
\begin{equation}
  \label{eq:gmam_descent}
  -\lambda H_{\vartheta\vartheta}D_\varphi\hat S(\varphi) = \lambda^2 \varphi'' 
  - \lambda H_{\vartheta\varphi}\varphi' +
  H_{\vartheta\vartheta}H_\varphi 
  + \lambda \lambda' \varphi'\,.
\end{equation}
This is derived as Proposition 3.1 in Appendix E of
\cite{heymann-vanden-eijnden:2008}, and we will show below how this
expression can be intuitively understood. Since
$H_{\vartheta \vartheta}$ is assumed to be positive definite and
$\lambda\ge0$, we can use~\eqref{eq:gmam_descent} directly to compute the
solution of~\eqref{eq:9} that minimizes $\hat S(\varphi)$ via a
relaxation method in virtual time $\tau$, that is, using the equation:
\begin{equation}
  \label{eq:10}
  \begin{aligned}
    \frac{\partial \varphi}{\partial\tau} & = - \lambda
    H_{\vartheta\vartheta}D_\varphi\hat S(\varphi)\\
    & = \lambda^2 \varphi'' - \lambda H_{\vartheta\varphi}\varphi' +
    H_{\vartheta\vartheta}H_\varphi + \lambda \lambda' \varphi'\,.
  \end{aligned}
\end{equation}
This equation is the main equation used in the original gMAM. Note that
the computation of the right hand-side of this equation requires the
computation of $H_\varphi$, $H_{\vartheta\varphi}$ and
$H_{\vartheta \vartheta}$, where the second derivatives of the
Hamiltonian possibly become unsightly for more complicated systems
that arise naturally when trying to use gMAM in practical
applications. In Sect.~\ref{sec:sgmam} we propose a simplification of
this algorithm that reduces the terms necessary to only first order
derivatives of the Hamiltonian, $H_\vartheta$ and $H_\varphi$.

Coming back to~\eqref{eq:gmam_descent}, it can be intuitively
understood by using the associated Hamiltonian system. Consider a
reparametrization of the original minimizer
$\varphi_*(s(t)) = \phi_*(t)$. In the following we are using a dot in
order to denote partial derivatives with respect to time and a prime
in order to denote a partial derivative with respect to the
parametrization $s$, hence $\dot v \equiv \partial v/\partial t$ and
$v' \equiv \partial v/\partial s$. With this notation, we find for
$\lambda^{-1} = t'(s)$ that $\dot\phi_* = \lambda\varphi_*'$ as well
as $\dot\phi_* = H_\theta, \dot\theta_* = -H_\phi$, and therefore
\begin{align*}
  \ddot\phi_* &= H_{\theta\phi} \dot\phi_* + H_{\theta\theta} \dot\theta_*\\
  &= \lambda H_{\theta\phi} \varphi_*' - H_{\theta\theta} H_\phi
\end{align*}
but also, since $\partial/\partial t = \lambda\, \partial/\partial
s$,
\begin{align*}
  \ddot\phi_* &= \partial(\lambda \varphi_*')/\partial t\\
  &= \lambda \lambda' \varphi_*' + \lambda^2 \varphi_*''
\end{align*}
so in total
\begin{equation*}
  - \lambda \lambda' \varphi_*' + \lambda H_{\theta\phi} \varphi_*' - 
  H_{\theta\theta} H_\phi - \lambda^2 \varphi_*'' = 0 = \lambda H_{\theta\theta}D_\varphi\hat S(\varphi)\,,
\end{equation*}
i.e. indeed the gradient vanishes at the minimizer.

\subsection{A simplified gMAM}
\label{sec:sgmam}

In contrast to the previous section, we start from the
form~\eqref{eq:gmam-2} of the geometric action. We want to solve the
mixed optimization problem\index{simplified geometric minimum action
  method}, i.e. find a trajectory $\varphi_*$ such that
\begin{equation}
  \varphi_* = \argmin_{\varphi\in\hat{\mathcal{C}}_{x,y}} \sup_{\vartheta: H(\varphi,\vartheta)=0} E(\varphi,\vartheta)\,,
\end{equation}
where
\begin{equation}
  E(\varphi,\vartheta) = \int_0^1 \dotprod{\varphi'}{\vartheta}\,ds\,.
\end{equation}
Let
\begin{equation}
  \label{eq:17}
  E_*(\varphi) =  \sup_{\vartheta:H(\varphi,\vartheta)=0} E(\varphi,\vartheta)
\end{equation}
and $\vartheta_*(\varphi)$ such that $E_*(\varphi) = E(\varphi,
\vartheta_*(\varphi))$.  This implies that $\vartheta_*$ fulfills the
Euler-Lagrange equation\index{Euler-Lagrange equation} associated with
the constrained optimization problem in~\eqref{eq:17}, that is,
\begin{equation}
  \label{eq:optim-theta} 
  D_\vartheta E(\varphi, \vartheta_*) = \mu H_\vartheta(\varphi,\vartheta_*)\,,
\end{equation}
where on the right-hand side $\mu(s)$ is the Lagrange
multiplier\index{Lagrange multiplier} added to enforce the constraint
$H(\varphi,\vartheta_*)=0$. In particular, at $\vartheta=\vartheta_*$,
we have
\begin{equation}
  \mu = \frac{\|D_\vartheta E\|^2}{\ddotprod{D_\vartheta E}{H_\vartheta}} = \frac{\|\varphi'\|^2}{\ddotprod{\varphi'}{H_\vartheta}}\,,
\end{equation}
where the inner product $\ddotprod{\cdot}{\cdot}$ and its induced norm
$\|\cdot\|$ can be chosen appropriately, for example as
$\dotprod{\cdot}{\cdot}$ or
$\dotprod{\cdot}{H_{\vartheta\vartheta}^{-1} \,\cdot}$.

At the minimizer $\varphi_*$, the variation of $E_*$ with respect to $\varphi$
vanishes. Using \eqref{eq:optim-theta} we conclude
\begin{align}
  0 = D_\varphi E_*(\varphi_*) &= D_\varphi E(\varphi_*,\vartheta_*) 
                           + \left[D_\vartheta E D_\varphi \vartheta
                           \right]_{(\varphi,\vartheta)=(\varphi_*,\vartheta_*)}\nonumber\\
                         &= -\vartheta_*' + \mu \left[H_\vartheta 
                           D_\varphi \vartheta\right]_{(\varphi,\vartheta)=(\varphi_*,\vartheta_*)}
                           \nonumber\\
                         &= -\vartheta_*' - \mu
                           H_\varphi(\varphi_*,\vartheta_*)\,,
                           \label{eq:optim-psi}
\end{align}
where in the last step we used $H(\varphi, \vartheta_*)=0$ and
therefore $$H_\varphi(\varphi, \vartheta_*) = -
H_\vartheta(\varphi,\vartheta_*) D_\varphi \vartheta.$$ Multiplying
the gradient \eqref{eq:optim-psi} with any positive definite matrix as
pre-conditioner\index{pre-conditioner} yields a descent direction. It
is necessary to choose $\mu^{-1}$ as pre-conditioner to ensure
convergence around critical points, where $\varphi'=0$.

Summarizing, we have reduced the minimization of the geometric action
into two separate tasks:
\begin{enumerate}
\item \index{conjugate momentum}For a given $\varphi$, find $\vartheta_*(\varphi)$ by solving the
  constrained optimization problem\index{constrained optimization problem}
  \begin{equation}
    \label{eq:theta-star}
    \vartheta_*(\varphi) = \argmax_{\vartheta, H(\varphi,\vartheta)=0} E(\varphi,\vartheta)\,,
  \end{equation}
  which is equivalent to solving
  \begin{equation}
    \label{eq:22}
    D_\vartheta E(\varphi,\vartheta_*) = \varphi' = \mu H_\vartheta(\varphi,\vartheta_*)
  \end{equation}
  for $(\mu,\vartheta_*)$ under the constraint $H(\varphi,\vartheta_*)=0$. This
  can be done via
  \begin{itemize}
  \item gradient descent;
  \item a second order algorithm for faster convergence
    (e.g. Newton-Raphson, as employed in
    \cite{heymann-vanden-eijnden:2008});
  \item in many cases, analytically (see below).
  \end{itemize}
\item Find $\varphi_*$ by solving the optimization problem
  \begin{equation}
    \varphi_* = \argmin_{\varphi\in\hat{\mathcal{C}}_{x,y}} E_*(\varphi)\,,
  \end{equation}
  for example by pre-conditioned gradient descent, using as direction
  \begin{equation}
    \label{eq:descent}
    -\mu^{-1} D_\varphi E_* = \mu^{-1}\vartheta_*'(\varphi) + H_\varphi(\varphi,\vartheta_*(\varphi))\,,
  \end{equation}
  with $\mu^{-1}$ as pre-conditioner. The constraint on the
  parametrization, e.g. ${|\varphi'|=\textrm{const}}$, must be fulfilled
  during this descent (see below).
\end{enumerate}

\subsection{Connection to gMAM}

The problem of finding $\vartheta_*(\varphi)$ is equivalent
to~\eqref{eq:find-theta-gmam} from gMAM and the same methods are
applicable. In particular note that the Lagrange multiplier $\mu$ which enforces
$H(\varphi_*,\vartheta_*)=0$ is identical to $\lambda^{-1}$.

It is also easy to see that, at $(\varphi_*,\vartheta_*)$, the combined optimization problem
$\{D_\vartheta E=\mu H_\vartheta, D_\varphi E_* = 0\}$ is identical to the geometric
equations of motion,
\begin{equation}
  \label{eq:eqnmotion-g}
  \begin{cases}
    D_\vartheta E = \varphi' = \mu H_\vartheta\\
    D_\varphi E_* = -\vartheta' - \mu H_\varphi = 0\,.
  \end{cases}
\end{equation}

On the other hand, none of the formulas in the above section use
higher derivatives of the Hamiltonian: Only $H_\varphi$ and $H_\vartheta$
are needed, which is a big simplification. This is obviously also true
for the equations of motion~\eqref{eq:eqnmotion} and their geometric
variant~\eqref{eq:eqnmotion-g}, which is the basis for the efficiency
of algorithms like \cite{chernykh-stepanov:2001,
  grafke-grauer-schaefer-etal:2014, grafke-grauer-schindel:2015}.

\subsection{Simplifications for SDEs with additive noise\index{SDE with additive noise}\index{additive noise}}
\label{ssec:additive}

For an SDE of the form
\begin{equation}
  dX = b(X)dt + \sqrt{\epsilon}\, dW\,,
\end{equation}
where $\sigma = \textrm{Id}$, the equations of gMAM become
significantly simpler. In the following, we derive explicit
expressions for this case, as it arises in numerous applications.

The corresponding Hamiltonian\index{Hamiltonian for additive noise} is given by
\begin{equation}
  \label{eq:H-quad-id}
  H(\varphi,\vartheta)=\dotprod{b}{\vartheta} + \frac12 \dotprod{\vartheta}{\vartheta}=0
\end{equation}
and we find directly
\begin{equation*}
  H_\varphi = (b_\varphi)^T \vartheta,\qquad H_\vartheta = b + \vartheta\,.
\end{equation*}
In many cases, we consider exits from stable fixed points of the
deterministic system where we have $H=0$ which, if we also use $D_\vartheta E
= \mu H_\vartheta$, permits to conclude that
\begin{equation}
  |H_\vartheta|^2=|b+\vartheta|^2 = |b|^2 + 2\dotprod{b}{\vartheta} +
  \dotprod{\vartheta}{\vartheta}
  =|b|^2+2H=|b|^2\,.
\end{equation}
As a result
\begin{equation}
  \label{eq:mu-sde}
  \mu = \frac{|D_\vartheta E|}{|H_\vartheta|} = \frac{|\varphi'|}{|b+\vartheta|}
  = \frac{|\varphi'|}{|b|}\,,
\end{equation}
i.e. we can compute $\mu$ without the knowledge of $\vartheta$. On the other
hand \eqref{eq:22} implies
\begin{equation}
  \label{eq:theta-sde}
  \varphi' = \mu H_\vartheta = \mu(b+\vartheta) \quad\Rightarrow\quad \vartheta = 
  \mu^{-1}\varphi' - b\,.
\end{equation}
The whole algorithm therefore reduces to the gradient descent
\begin{equation} \label{eq:gradient_descent}
  \frac{\partial \varphi}{\partial \tau} = \mu^{-1}\vartheta_*' + (b_\varphi)^T \vartheta_*\,,
\end{equation}
with $\mu, \vartheta_*$ given by~\eqref{eq:mu-sde}
and~\eqref{eq:theta-sde}. Examples in this class will be
treated in Secs.~\ref{ssec:maierstein}, \ref{ssec:acch-reduced}, and
\ref{ssec:meta-stable-climate} below.

\subsection{Simplifications for general SDEs (multiplicative noise)}
\label{ssec:multiplicative}

As a slightly more complicated case, consider the following SDE with
multiplicative noise\index{SDE with multiplicative
  noise}\index{multiplicative noise}:
\begin{equation}
  dX= b(X) \,dt + \sqrt{\epsilon} \sigma(X) \,dW\,,
\end{equation}
where $a(\varphi) = \sigma(\varphi) \sigma^\dagger(\varphi)$. Then the
Hamiltonian\index{Hamiltonian for multiplicative noise} reads
\begin{equation}
  H(\varphi,\vartheta) = \dotprod{b}{\vartheta} + \tfrac12 \dotprod{\vartheta}{a \vartheta}
\end{equation}
and
\begin{equation}
  H_\varphi = (b_\varphi)^T\vartheta + \tfrac12 \dotprod{\vartheta}{(a_\varphi)\vartheta},\qquad H_\vartheta = b + a \vartheta\,.
\end{equation}
Defining an inner product and norm induced by the correlation,
$\dotprod{u}{v}_a=\dotprod{u}{a^{-1}v}$ and $|u|_a =
\dotprod{u}{u}_a^{1/2}$ yields, as before,
\begin{equation}
  \label{eq:mu-mul}
  |H_\vartheta|_a = |b|_a \qquad \Rightarrow \qquad \mu =
  \frac{|\varphi'|_a}{|b|_a}
\end{equation}
and
\begin{equation}
  \label{eq:theta-mul}
  \vartheta = a^{-1}(\mu^{-1}\varphi' - b).
\end{equation}
In the case of multiplicative noise, the algorithm therefore reads
\begin{equation}
  \frac{\partial \varphi}{\partial \tau} = \mu^{-1}\vartheta_*' + \left((b_\varphi)^T \vartheta_* +\tfrac12 \dotprod{\vartheta_*}{(a_\varphi) \vartheta_*}\right)\,,
\end{equation}
with $\mu, \vartheta_*$ given by~\eqref{eq:mu-mul}
and~\eqref{eq:theta-mul}. An example in this will be treated in
Sect.~\ref{ssec:voter-model}.

It is also worth pointing out that we encounter difficulties as soon
as the noise correlation $a$ is not invertible. This is equivalent to
stating that some degrees of freedom are not subject to noise and thus
behave deterministically. The adjoint field $\vartheta$ has to be equal
to zero on these modes, and they fulfill the deterministic equation
$\varphi' = b$ exactly. This translates into additional constraints
for the minimization procedure, which have to be enforced numerically.

\subsection{Comments on improving the numerical efficiency}

To increase the numerical efficiency\index{numerical
  efficiency} of the algorithm, some alterations are possible:
\begin{itemize}
\item Arc-length parametrization\index{arc-length parametrization},
  $|\varphi'|=\textrm{const}$, can be enforced trivially and without
  introducing a stiff Lagrange multiplier term by interpolation along
  the trajectory every (or every few) iterations. As additional
  benefit of this method all terms of the relaxation dynamics which
  are proportional to $\varphi'$ can be discarded, as they are
  canceled by the reparametrization. This is of particular use in
  applications that involve PDEs (see Sect.~\ref{sec:pde}), as shown
  in examples below.
\item Stability\index{numerical stability} in the relaxation parameter
  can be greatly increased if one treats the stiffest term of the
  relaxation equation implicitly. In ODE systems, the stiffest term
  usually is $H_{\vartheta\vartheta}^{-1} \varphi''$, which is
  contained in $\vartheta'$. For simplicity of implementation, it is
  sufficient to compute $\vartheta_*$ in the usual way, apply
  $\vartheta_*'$ in the descent step, but subtract
  $H_{\vartheta\vartheta}^{-1} \varphi_n''$ and add
  $H_{\vartheta\vartheta}^{-1} \varphi_{n+1}''$ here. This also works
  in the case of general Hamiltonians, where the dependence of
  $\vartheta_*$ on $\varphi'$ is less obvious.

  In our implementation, the relaxation step is conducted by computing
  \begin{equation}
    \label{eq:strang}
    \varphi_{n+1} = \left(1 - h \mu^{-2} H_{\vartheta\vartheta}^{-1}
    \partial_s^2\right)^{-1} R_n\,,
  \end{equation}
where
$$
R_n = \left(\varphi_n + h (\mu^{-1}\vartheta_*'(\varphi_n) + 
    H_\varphi(\varphi_n,\vartheta_*(\varphi_n)) - \mu^{-2}H_{\vartheta\vartheta}^{-1} \varphi_n'')
    \right)\,.
$$
  This division into an implicit treatment of the stiffest term and
  explicit treatment of the rest is the simplest case of
  Strang splitting \cite{strang:1968} and the implementation of
  \eqref{eq:strang} is only first order accurate. The splitting can be
  taken to arbitrary order \cite{yoshida:1990} under additional
  computational cost.

  Note that the above modification, while increasing efficiency, at
  the same time increases complexity, as the computation of the second
  derivative $H_{\vartheta\vartheta}$ becomes necessary. In practice, if the
  Hamiltonian is not too complex, we find that the benefits outweigh
  the implementation costs, and some problems, especially PDE systems,
  are not tractable at all with the inefficient but simpler choice of
  explicit relaxation. If the PDE system contains higher-order spatial
  derivatives, even more terms should possibly be treated with a
  stable integrator, as is discussed in the next section.
\item Depending on the problem, it might be beneficial to choose a
  different scalar product in the descent. In case of traditional
  gMAM, the descent is done using
  $\langle \cdot, (\mu^2 H_{\vartheta\vartheta})^{-1} \,\cdot \rangle$, but
  other choices are also feasible. Note that it is possible to choose
  the metric such that at least one term at the right-hand side
  disappears, as it becomes parallel to the trajectory and is canceled
  by reparametrization, as outlined above.
\item Some insight about the nature of the transition can be obtained
  by first finding the heteroclinic orbits\index{heteroclinic orbit}
  defined geometrically as
  \begin{equation}
    \varphi' \parallel b(\varphi)\,.
  \end{equation}
  This calculation can be done very efficiently even for complicated
  problems via the string method
  \cite{e-ren-vanden-eijnden:2002}. Even though the heteroclinic orbit
  differs from the transition path for systems that violate detailed
  balance, it \emph{does} correctly predict the transition from the
  saddle point onward (the ``downhill'' portion, which happens
  deterministically). The method put forward here can then be used to
  find the transition path up to the saddle (the ``uphill'' portion)
  only. If there are several saddles to be taken into account, it is
  not known \emph{a priori} which one will be visited by the
  transition pathway. In this case, the strategy has to be modified
  accordingly, for example by computing one heteroclinic orbit per
  saddle. To highlight the relation between the string and the
  minimizer, we compute and compare the two in many of the
  applications below. We denote with ``string'' the heteroclinic
  orbits connecting the fixed points to the saddle point of relevance
  found via the string method.
\end{itemize}

\subsection{SPDEs with additive noise}
\label{sec:pde}

\index{SPDEs with additive noise}\index{additive noise}In this
section, we discuss the application to SPDE systems. For simplicity,
we focus on the case of SPDEs with additive noise that can be
written formally as
\begin{equation} \label{eq:basic_spde}
  U_t = B(U) + \sqrt{\epsilon}\, \eta(x,t)\,,
\end{equation}
where the drift term is given by the operator $B(U)$ and $\eta$
denotes spatio-temporal white-noise. It is a non-trivial task to make
mathematical sense of such SPDEs under spatially irregular noise due
to the possible ill-posedness of non-linear terms, especially if the
spatial dimension is higher than one. This may require to renormalize
the equation, which can be done rigorously in certain cases using the
theory of regularity structures \cite{hairer:2014}. The
renormalization procedure typically involves mollifying the noise term
on a scale $\delta$, and adding terms in the equation that
counterbalance divergences that may occur as one lets $\delta\to0$. In
the context of LDT, the main issue is whether these renormalizing
terms subsist if we also let
$\epsilon\to0$. In~\cite{hairer-weber:2015}, it was shown in the
context of the stochastic Allen-Cahn equation in 2 or 3 spatial
dimensions that the action of the mollified equation converges towards
the action associated with the (possibly formal) equation
in~\eqref{eq:basic_spde} in which the noise is white-in-space provided
that $\epsilon$ is sent to zero fast enough as $\delta\to0$. This
action reads
\begin{equation}
  \label{eq:7}
  S_T(\phi)  = \frac12\int_0^T \|\phi_t - B(\phi)\|_{L^2}^2dt\,,
\end{equation}
where $\|\cdot\|_{L^2}$ denotes the $L^2$-norm. This leads to
expressions for the geometric action that are similar to those
in~\eqref{eq:gmam} but with the Euclidean inner product replaced by
the $L^2$-inner product. In the sequel we will not dwell further on
these mathematical issues and always assume that \eqref{eq:7} and the
associated geometric action are the relevant one to study.

The gradient descent for the minimizer of this geometric action is
similar to the one in (\ref{eq:gradient_descent}) but with the term
$(b_\varphi)^T$ replaced by the functional derivative of the operator
$B$ with respect to $\varphi$.
\begin{equation} \label{eq:gradient_descent_pde}
\frac{\partial \varphi}{\partial \tau} = \mu^{-1}\vartheta_*' 
+ \left(D_\varphi B \right)^T \vartheta_*\,.
\end{equation}
In practice, however, this equation needs to be rewritten in order to
allow for numerical stability. This is due to the fact that the scheme
will contain derivatives of high orders, and their corresponding
stability condition (CFL condition) will limit the rate of convergence
of the scheme. We therefore want to treat the most restrictive terms
either implicitly or with exponential integrators. To this end, let us
focus on the following class of problems where the drift $B$ can be
written as
\begin{equation}
  \label{eq:split}
  B = L\varphi + R(\varphi)\,,
\end{equation}
where $L$ is a linear self-adjoint operator containing higher-order
derivatives that does not depend on time explicitly, and $R(\varphi)$ is
the rest, possibly nonlinear. Recall that $\vartheta_*$ can be computed
from $\varphi'$ via
\begin{equation}
  \vartheta_* = \mu^{-1}\varphi' - B = \mu^{-1}\varphi'- L\varphi - R(\varphi)\,.
\end{equation}
On the other hand, we have also a term proportional to $L$ in
\begin{equation}
  D_\varphi B =D_\varphi R + L
\end{equation}
and, therefore, the relaxation formula (\ref{eq:gradient_descent_pde})
for $\varphi$ actually contains a term $L^2 \varphi$. If $L$ contains
higher-order derivatives, this term will likely be the most
restrictive in terms of numerical stability. It is therefore
advantageous to treat it separately. Introducing an auxiliary variable
$\tilde \vartheta_*$ defined by
\begin{equation}
  \tilde \vartheta_* = \mu^{-1}\varphi' - R(\varphi) = \vartheta_* + L \varphi
\end{equation}
we can rewrite the relaxation formula as
\begin{eqnarray*} 
  \frac{\partial \varphi}{\partial \tau} &=& \mu^{-1}\vartheta_*' + \left(D_\varphi B \right)^T \vartheta_*\, \\
 &=& \mu^{-1}\tilde \vartheta_*' - \mu^{-1}L \varphi' + \left(D_\varphi R \right)^T \vartheta_* + L \vartheta_* \\
 &=& \tilde \mu^{-1}\vartheta_*' - \mu^{-1}L \varphi' + \left(D_\varphi R \right)^T \vartheta_* + L (\tilde \vartheta_* -  L \varphi) \\
 &=& \mu^{-1}\tilde \vartheta_*' - \mu^{-1}L \varphi' + \left(D_\varphi R \right)^T \vartheta_* + L \tilde \vartheta_* - L^2 \varphi \\
 &=& \mu^{-1}\tilde \vartheta_*' + \left(D_\varphi R \right)^T \vartheta_* - L R(\varphi) - L^2 \varphi\,.
\end{eqnarray*}
The term $L^2 \varphi$ is now separated and can be treated
independently. Since it is linear by definition, it can be treated
very efficiently with an integrating factor by employing exponential
time differencing\index{exponential time differencing} (ETD)
\cite{beylkin-keiser-vozovoi:1998}. For an equation with a
deterministic term of the form \eqref{eq:split}, multiplying by the
integrating factor $e^{-L\tau}$ and integrating from $\tau_n$ to
$\tau_{n+1}=\tau_n+h$, one obtains the \emph{exact} formula
\begin{equation}
  \varphi_{n+1} = e^{Lh}\varphi_n + e^{Lh} \int_0^h e^{-L\tau} R(\varphi(t_n+\tau))\,d\tau\,,
\end{equation}
which can be approximated by
\begin{equation}
  \label{eq:ETD}
  \varphi_{n+1} = e^{Lh}\varphi_n + (e^{Lh}-\textrm{Id})L^{-1} R(\varphi_n)\,,
\end{equation}
when treating the linear part of the equation exactly and
approximating the integral to first order. This scheme can be taken to
higher order \cite{cox-matthews:2002} and its stability improved
\cite{kassam-trefethen:2005}, but a first order scheme proved to be
sufficient for the examples given below. For the descent
\eqref{eq:gradient_descent_pde} we want to treat the stiffest part
$-L^2 \varphi$ with ETD, so the integrating factor here becomes
$e^{-L^2 \tau}$.

A complete relaxation step then consists of
\begin{enumerate}
\item compute $\vartheta_*$ and $\tilde \vartheta_*$ using the explicit
  formulas
  \begin{equation*}
    \tilde \vartheta_* = \mu^{-1} \varphi' - R(\varphi) , \qquad \vartheta_* = \tilde
    \vartheta_* - L \varphi\,;
  \end{equation*}
\item compute the explicit step
  \begin{equation*}
    \xi = \mu^{-1}\tilde \vartheta_*' + \left(D_\varphi R \right)^T \vartheta_* -L R(\varphi)  - \mu^{-2}H_{\vartheta\vartheta}^{-1} \varphi_n''\,,
  \end{equation*}
  where as in the SDE case, if needed, we can subtract the term
  $\mu^{-2}H_{\vartheta\vartheta}^{-1} \varphi_n''$ to treat it implicitly
  later;
\item perform an ETD step
  \begin{equation*}
    \bar \varphi = e^{-L^2 h} \varphi_n - (e^{-L^2 h}-\textrm{Id})(L^2)^{-1} \xi\,;
  \end{equation*}
\item apply the second derivative in arc-length direction implicitly,
  \begin{equation*}
    \varphi_{n+1} = (1-h\mu^{-2}H_{\vartheta\vartheta}^{-1}\partial_t^2)^{-1} \bar \varphi\,.
  \end{equation*}
  \label{itm:implicit}
\end{enumerate}

Note that the integral factors $e^{-L^2 h}$ and
$(e^{-L^2h}-\textrm{Id})(\mu L^2)^{-1}$ are possibly costly to
compute, as they contain matrix-exponentials and inversions. However,
the computation can be done once before starting the iteration, so
that the associated computational cost becomes negligible. In
contrast, this is not true in general for the implicit step
\ref{itm:implicit}, since $\mu^{-2}H_{\vartheta\vartheta}^{-1}$ might depend
on the fields in a complicated way and has to be recomputed at every
iteration.

\section{Illustrative applications}
\label{sec:applications}

In what follows we apply our simplified gMAM to the series of examples
listed in the introduction. These examples illustrate specific
questions encountered in practical applications arising in a variety
of fields, in which the computation of the rate and mechanism of
transitions is of interest. Note that all these examples involve
non-equilibrium systems whose dynamics break detailed balance, so that
simpler methods of computation are not readily available.

In the following, we will break our notation convention and instead
use the notation of the respective fields to minimize confusion.

\subsection{Maier-Stein model}
\label{ssec:maierstein}
\index{Maier-Stein model}

\begin{figure}[tb]
  \begin{center}
    \includegraphics[width=0.48\textwidth]{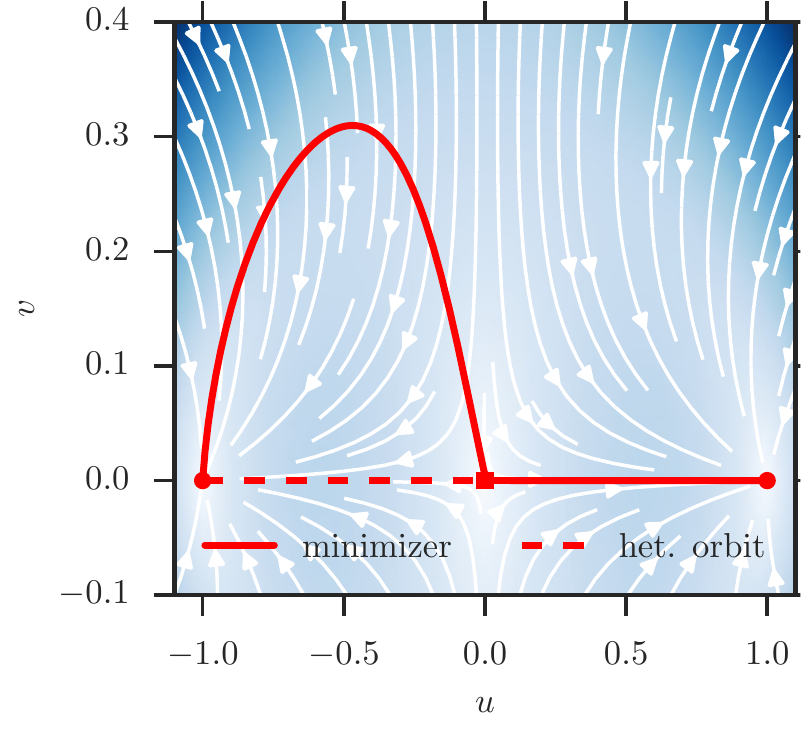}\hfill
    \includegraphics[width=0.48\textwidth]{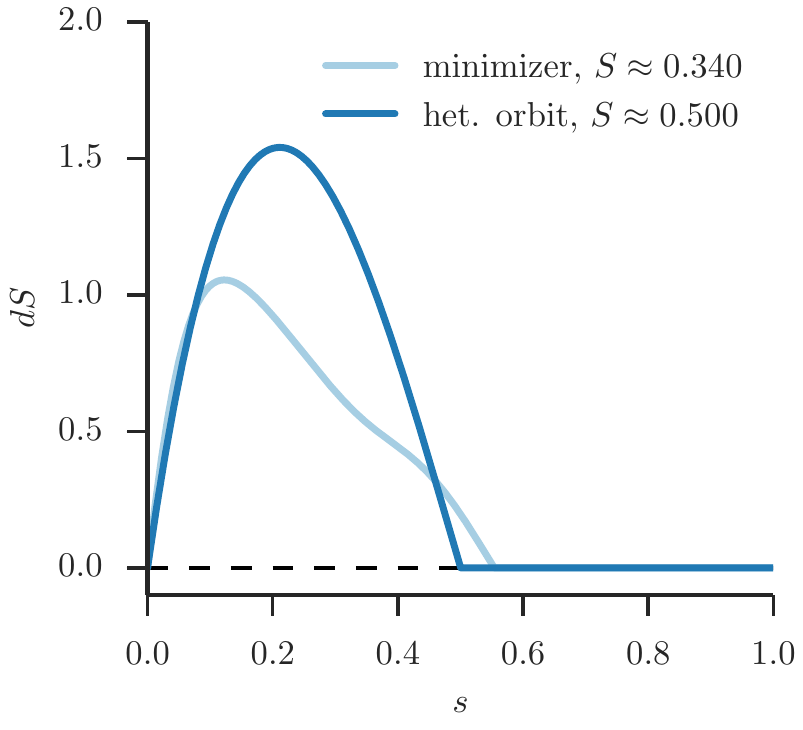}
  \end{center}
  \caption{Maier-Stein model, $\beta=10$. Left: PML and heteroclinic
    orbit. The arrows denote the direction of the deterministic flow,
    the shading its magnitude. The solid line depicts the minimizer,
    the dashed line the heteroclinic orbit. Dots are located at the
    fixed points (circle: stable; square: saddle). Right: Action
    density along the minimizer and the heteroclinic
    orbit.\label{fig:maierstein}}
\end{figure}

Maier and Stein's model \cite{maier-stein:1996} is a simple system
often used as benchmark in LDT calculations. It reads
\begin{equation}
  \label{eq:maierstein}
  \begin{cases}
    du = (u-u^3-\beta u v^2) dt + \sqrt{\epsilon}dW_u\\
    dv = -(1+u^2)v dt + \sqrt{\epsilon}dW_v\,,
  \end{cases}
\end{equation}
where $\beta$ is a parameter. For all values of $\beta$, the
deterministic system has the two stable fixed points,
$\varphi_-=(-1,0)$ and $\varphi_+=(1,0)$, and a unique unstable
critical point $\varphi_s=(0,0)$. However it satisfies detailed
balance only for $\beta=1$. In this case, we can write the drift in
gradient\index{gradient system} form, $b(\varphi) = \nabla_\varphi
U(\varphi)$, and the minimizers of the geometric action that connects
$\varphi_- $ to $\varphi_+$ and \textit{vice-versa} are the
time-reverse of each other and lie on the location of the heteroclinic
orbit where $\varphi' \parallel \nabla U$.  Here, we use $\beta=10$,
in which case detailed balance is broken and the forward and backward
transition pathways are no longer identical. Since the noise is
additive, the system~\eqref{eq:maierstein} falls into the category
discussed in Sect.~\ref{ssec:additive} (i.e. additive noise) and can
be solved with the simplest variant of the algorithm. The minimizer of
the action connecting $\varphi_- $ to $\varphi_+$ and the value of the
action along it are shown in Fig.~\ref{fig:maierstein}. Since the
system is invariant under the transformation $v\to-v$, there is also a
minimizer with identical action in the $v<0$ half-plane. Similarly,
the paths from $\varphi_+$ to $\varphi_-$ can be obtained via the
transformation $u\to-u$. The numerical parameters used in these
calculations were $h=10^{-1}$, $N_s=2^{10}$, where $N_s$ denotes the
number of configurations along the transition trajectory or the number
of \emph{images}.

\subsection{Allen-Cahn/Cahn-Hilliard system}
\label{ssec:acch-reduced}
\index{Allen-Cahn equation}
\index{Cahn-Hilliard equation}

Pattern formation\index{pattern formation} in motile micro-organisms
is often driven by non-equilibrium forces, leading to visible patterns
in cellular colonies \cite{shapiro:1995,
  cates-marenduzzo-pagonabarraga-etal:2010}. For example,
\emph{E. coli} in a uniform suspension separates\index{phase
  separation} into a bacteria-rich and a bacteria-poor phase if the
swim speed decreases sufficiently rapidly with density
\cite{tailleur-cates:2008}. Here we study a model inspired by these
phenomena. We note that this model does not permit the thermodynamic
mapping used in \cite{tailleur-cates:2008}, so that understanding the
non-equilibrium transitions in the model requires minimization of the
geometric action of LDT.

\subsubsection{Reduced Allen-Cahn/Cahn-Hilliard model}

\begin{figure}[t]
  \includegraphics[width=200pt]{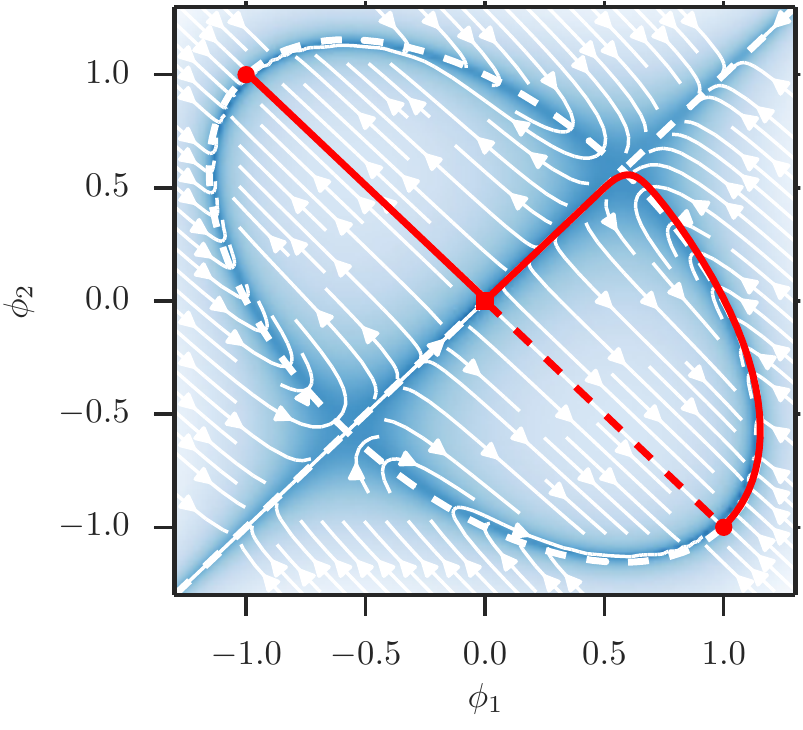}
  \caption{Allen-Cahn/Cahn-Hilliard toy ODE model, $\alpha=0.01$. The
    arrows denote the direction of the deterministic flow, the color
    its magnitude. The white dashed line corresponds to the slow
    manifold. The solid line depicts the minimizer, the dashed line
    the heteroclinic orbit. Markers are located at the fixed points
    (circle: stable; square: saddle).\label{fig:acchtoy-dynamics}}
\end{figure}

\begin{figure}
  \begin{center}
    \includegraphics[height=133pt]{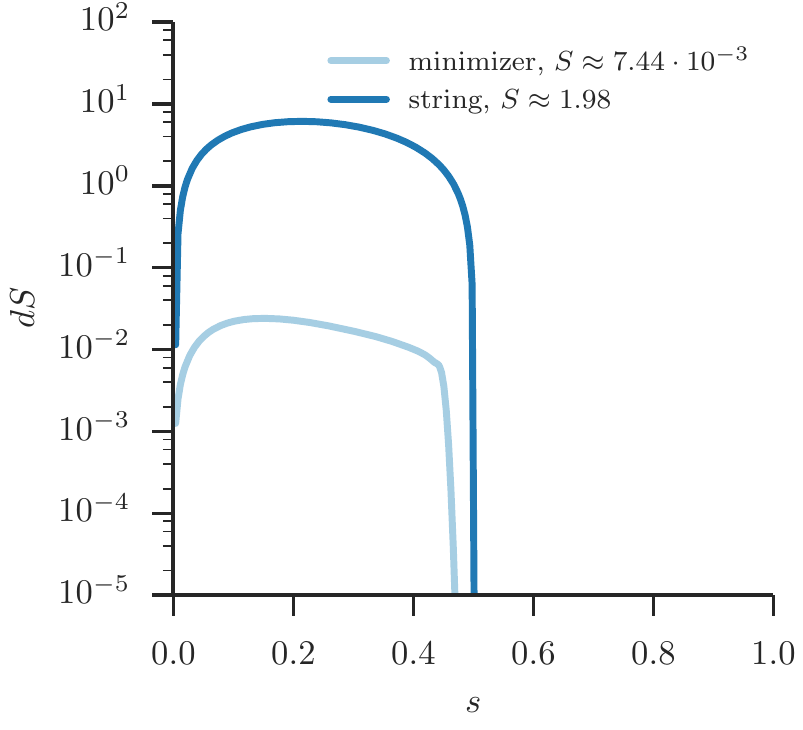}
    \includegraphics[height=133pt]{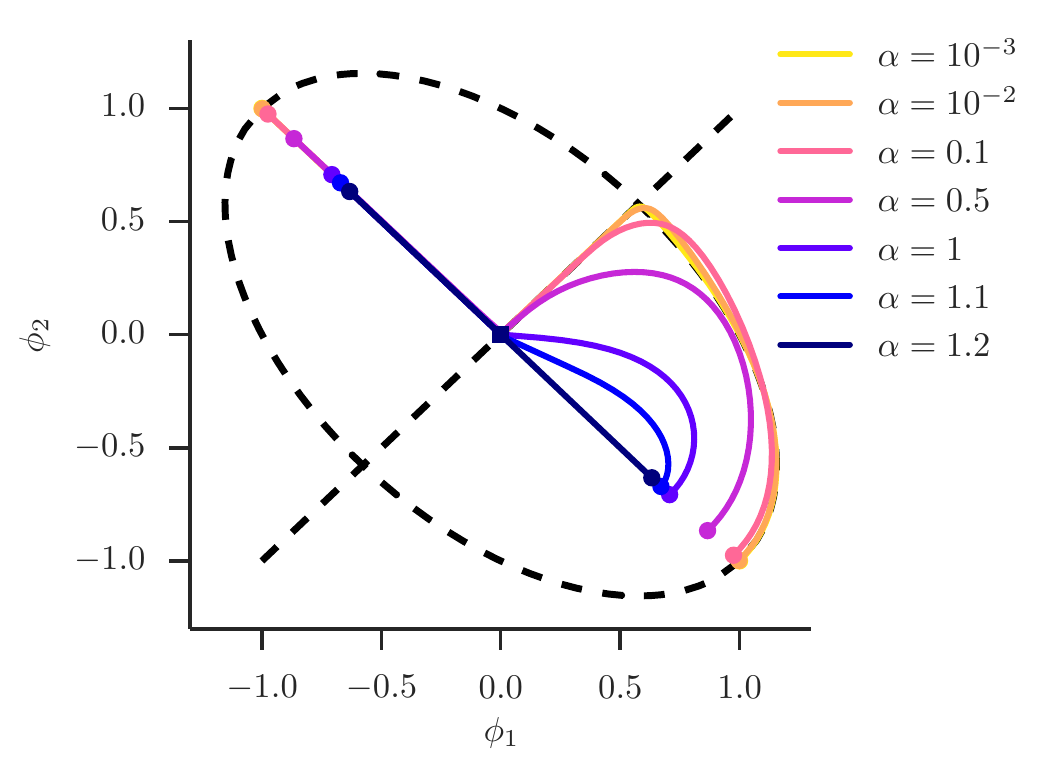}
  \end{center}
  \caption{Left: Action density along the path for the 2-dimensional
    reduced model. Path parameter is normalized to $s\in(0,1)$. For
    the second half of the transition, the action density is
    zero. Right: Minimizers of the action functional for different
    values of $\alpha$. For $\alpha\to0$, the minimizer approaches the
    slow manifold. Note that the switch to a straight line minimizer
    happens at a finite value $\alpha \approx
    1.12$.\label{fig:acchtoy-action}}
\end{figure}

Consider the SDE system
\begin{equation}
  \label{eq:acch-toy}
  d\phi= (\frac1\alpha Q (\phi-\phi^3) - \phi)dt + \sqrt{\epsilon} dW
\end{equation}
with $\phi=(\phi_1, \phi_2)$ and the matrix $Q=((1,-1),(-1,1))$. This
system does not satisfy detailed balance, as its drift is made of two
gradient terms with incompatible mobility operators (namely $Q$ and
$\textrm{Id}$). Model \eqref{eq:acch-toy} can be seen as a
2-dimensional reduction to a discretized version of the continuous
Allen-Cahn/Cahn-Hilliard model discussed later in Sect.~\ref{sec:acch}.

The deterministic flowlines of~\eqref{eq:acch-toy} are depicted in
Fig.~\ref{fig:acchtoy-dynamics}.  The deterministic dynamics has two
stable fixed points, $\phi_A=(-1,1)$ and $\phi_B=(1,-1)$, and an
unstable critical point, $\phi_S=(0,0)$, lying on the
separatrix\index{separatrix} where $\phi_1=\phi_2$ between the basins
of $\phi_A$ and $\phi_B$. The location of the heteroclinic orbits
connecting $\phi_S$ to $\phi_A$ and $\phi_B$ is a straight line
between these points. When $\alpha$ is small in~\eqref{eq:acch-toy},
there exists a ``slow manifold''\index{slow manifold}, comprised of
all points where $Q(\phi-\phi^3)=0$ which is shown as a white dashed
line in Fig.~\ref{fig:acchtoy-dynamics} . On this manifold, the
deterministic dynamics are of order $O(1)$, which is small in
comparison to the dynamics of the $Q$-term, which are of order
$O(1/\alpha)$. This suggests that for small enough $\alpha$ the
transition trajectory will follow this slow manifold on which the
drift is small, rather than the heteroclinic orbit, to escape the
basin of the stable fixed points. This is confirmed in
Fig.~\ref{fig:acchtoy-dynamics} where we show the action minimizer
connecting $\phi_B$ to $\phi_A$. As can be seen, the minimizer first
tracks the slow manifold, and it approaches the separatrix at a point
far from~$\phi_S$. It then follows closely the separatrix
towards~$\phi_S$ (which has to be part of the transition) to cross
into the other basin and then relax (deterministically)
towards~$\phi_A$.

The action along the minimizer and the paths made of the heteroclinic
orbits are depicted in Fig.~\ref{fig:acchtoy-action} (left). Notably,
due to its movement along the slow manifold, the action along the
minimizer is smaller by a factor of order $\alpha$. Minimizers for
different values of $\alpha$ are shown in
Fig.~\ref{fig:acchtoy-action} (right). Note that in the opposite limit
$\alpha\gg0$ the switch to a straight line happens at a finite value
$\alpha \approx 1.12$.

In these computations, we used $N_s = 2^{14}$, $h=10^{-2}$.

\subsubsection{Full Allen-Cahn/Cahn-Hilliard model}
\label{sec:acch}

Consider next the SPDE
\begin{equation}
  \label{eq:acch-model}
  \phi_t = \frac1\alpha P(\kappa \phi_{xx} + \phi - \phi^3) - \phi
  + \sqrt{\epsilon}\eta(x,t)\,,
\end{equation}
where $P$ is an operator with zero spatial mean and $\eta(x,t)$ a
spatio-temporal white-noise. This model is again of the form of two
competing gradient flows\index{competing gradient flows} with
different mobilities:
\begin{equation}
  \phi_t = -M_1 D_\phi V_1(\phi) -M_2 D_\phi V_2(\phi) + \sqrt{\epsilon}M_2^{1/2}\eta(x,t)\,,
\end{equation}
with 
\begin{subequations}
  \label{eq:form}
  \begin{align}
    V_1(\phi)&=\frac12 \kappa |\phi_x|^2 + \frac12 |\phi|^2 - \frac14
    |\phi|^4,& M_1&=\frac1\alpha P\\
    V_2(\phi) &= -\frac12 |\phi|^2,& M_2&=\textrm{Id}\,.
  \end{align}
\end{subequations}
For $P = -\partial_{x}^2$ the system is a mixture of a stochastic
Allen-Cahn \cite{allen-cahn:1972} and Cahn-Hilliard
\cite{cahn-hilliard:1958} equation. Here we will consider
$P(\phi) = \phi - \dashint\phi\, dx$, which is similar in most aspects
discussed below but simpler to handle numerically. We are again
interested in situations where $\alpha$ is small, and the time scales
associated with $V_1$ and $V_2$ differ significantly. In this case it
will turn out that transition pathways are very different
from the heteroclinic orbits, in that the separatrix between the
basins of attraction is approached far from the unstable critical
point of the deterministic system. This behavior is reminiscent of the
2-dimensional example discussed above, but in an SPDE setting.

\begin{figure}[tb]
  \begin{center}
    \includegraphics[width=215pt]{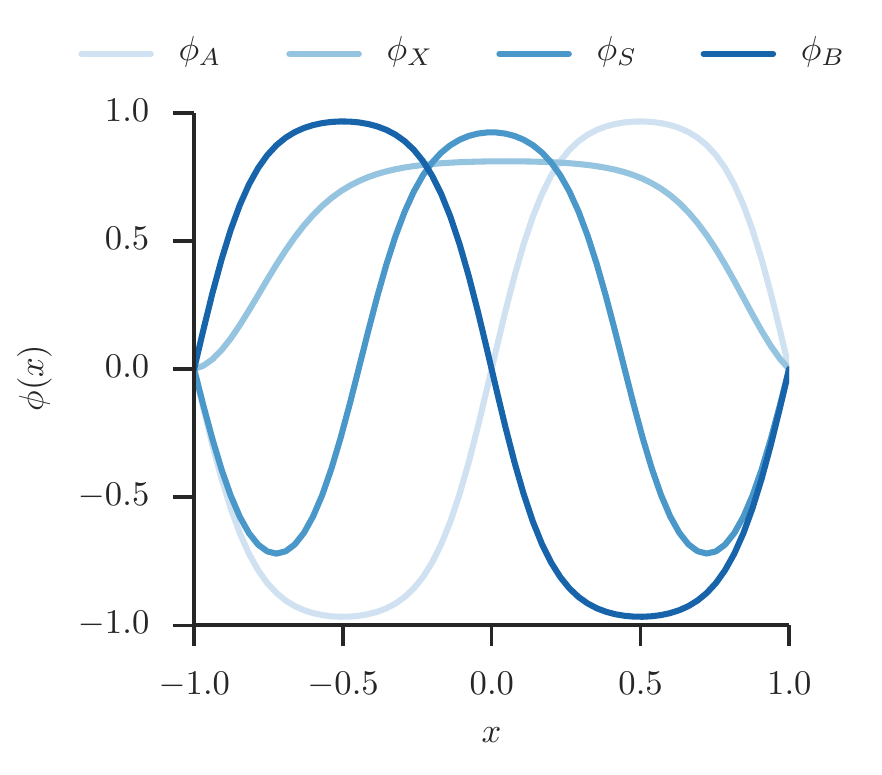}
  \end{center}
  \caption{The configurations $A, B, S, X$ in space: $\phi_A$ and
    $\phi_B$ are the two stable fixed points, $\phi_S$ is the unstable
    fixed point on the separatrix in between. At point $\phi_X$, the
    slow manifold intersects the separatrix. \label{fig:corners}}
\end{figure}

The fixed points of the deterministic ($\epsilon=0$) dynamics of
system ~\eqref{eq:acch-model} are the solutions of
\begin{equation}
  \label{eq:fixedp}
   P(\kappa\phi_{xx}+\phi-\phi^3)-\alpha \phi = 0\,.
\end{equation}
The only constant solution of this equation is the trivial fixed point
$\phi(x)=0$, whose stability depends on $\alpha$ and $\kappa$. In the
following, we choose $\alpha=10^{-2}$ and $\kappa = {2\cdot 10^{-2}}$,
in which case $\phi(x)=0$ is unstable.  The two stable fixed points
obtained by solving \eqref{eq:fixedp} for these values of $\alpha$ and
$\kappa$ are depicted in Fig.~\ref{fig:corners} as $\phi_A$ and
$\phi_B$, with $\phi_A = -\phi_B$. An unstable fixed point
configuration on the separatrix\index{separatrix} between $\phi_A$ and
$\phi_B$ is also shown as $\phi_S$.

\begin{figure}[tb]
  \begin{center}
    \includegraphics[width=0.48\textwidth]{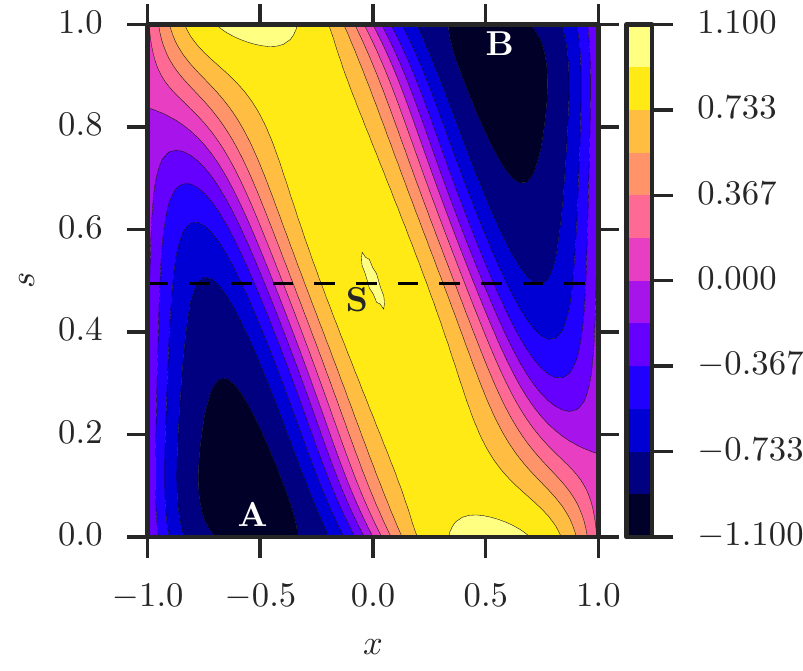}\hfill
    \includegraphics[width=0.48\textwidth]{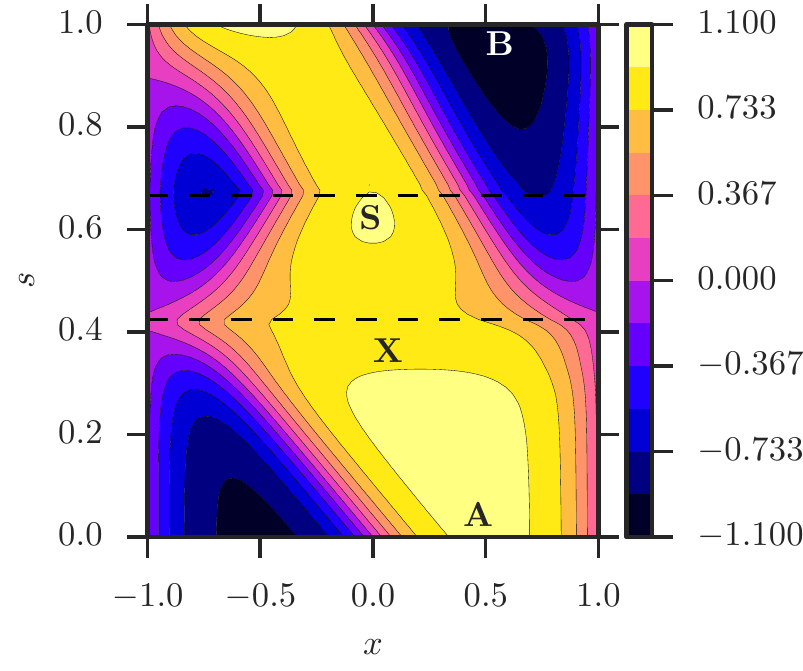}
  \end{center}
  \caption{Transition pathways between two stable fixed points of
    equation \eqref{eq:acch-model} in the limit $\epsilon\to0$. Left:
    heteroclinic orbit, defining the deterministic relaxation dynamics
    from the unstable point $S$ down to either $A$ or $B$. Right:
    Minimizer of the geometric action, defining the most probable
    transition pathway from $A$ to $B$, following the slow manifold up
    to $X$, where it starts to nearly deterministically travel close
    to the separatrix into $S$.\label{fig:transition}}
\end{figure}

For finite but small $\alpha$, the deterministic part
of~\eqref{eq:acch-model} has a ``slow manifold''\index{slow manifold}
made of the solutions of
\begin{equation}
  \label{eq:slowM}
  P(\kappa \phi_{xx}+\phi-\phi^3)=0\,.
\end{equation}
On this manifold the motion is driven solely by changing the mean via
the slow terms, $- \phi + \sqrt{\epsilon}\,\eta(x,t) $, on a
time-scale of order $O(1)$ in $\alpha$. After two integrations in
space, \eqref{eq:slowM} can be written as
\begin{equation}
  \kappa \phi_{xx} + \phi - \phi^3 = \lambda\,,
\end{equation}
where $\lambda$ is a parameter. As a result the slow manifold can be
described as one-parameter families of solutions parametrized by
$\lambda\in\mathbb{R}$ -- in general there is more than one family
because the manifold can have different branches corresponding to
solutions of~\eqref{eq:fixedp} with a different number of domain
walls. The configuration labeled as $\phi_X$ in Fig.~\ref{fig:corners}
shows the field at the intersection of one of these branches with the
separatrix.  Since the deterministic drift along the slow manifold is
small compared to the $O(1/\alpha)$ drift induced by the Cahn-Hilliard
term, one expects that the most probable transition pathway will use
this manifold as channel to escape the basin of attraction of the
stable fixed points $\phi_A$ or $\phi_B$.  This intuition is confirmed
by the numerics, as shown next.

\begin{figure}[tb]
  \begin{center}
    \includegraphics[width=215pt]{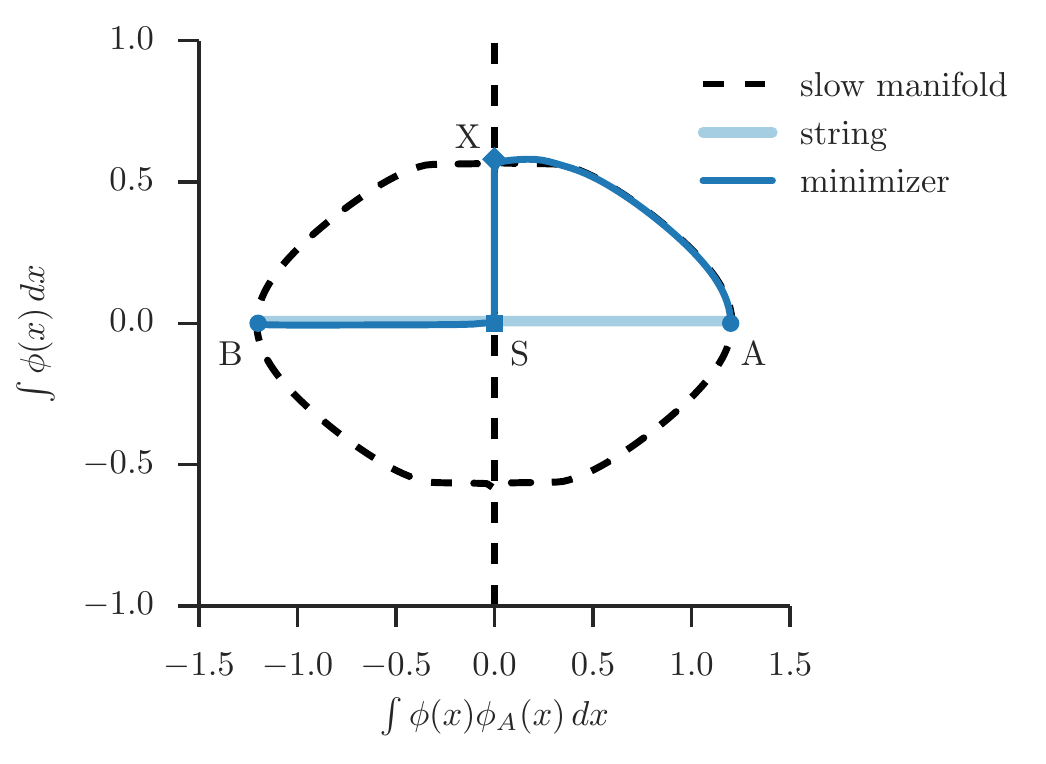}
  \end{center}
  \caption{Projection of the heteroclinic orbit and the minimizer of
    the action functional into a 2-dimensional plane. The
    $x$-direction is proportional to its component in the direction of
    the initial condition $\phi_A$ while the $y$-direction corresponds
    to its spatial mean. The stable fixed points are located at $A$
    and $B$, the unstable fixed point at $S$. The separatrix is the
    straight line $\int \phi(x)\phi_A(x)\,dx=0$. The heteroclinic
    orbit (light) travels $A\to S\to B$ in a horizontal line with
    vanishing mean, while the minimizer (dark) travels first along the
    slow manifold (dashed) $A\to X$ and then tracks the separatrix
    from $X$ to $S$. \label{fig:projection}}
\end{figure}
Fig.~\ref{fig:transition} (left) shows the heteroclinic orbit
connecting the two stable fixed points $\phi_A$ and $\phi_B$ to the
unstable configuration $\phi_S$. The mean is preserved along this
orbit, which involves a nucleation event at the boundaries followed by
domain wall motion through the domain. The unstable fixed point
$\phi_s$, denoted by $S$, which also demarcates the position at which the
separatrix is crossed, is the spatially symmetric configuration with a
positive central region and two negative regions at the
boundary. Locations $A$ and $B$ label the two stable fixed points
$\phi_A$ and $\phi_B$.

In contrast, Fig.~\ref{fig:transition} (right) shows the minimizer of
the geometric action, which is the most probable transition path as
$\epsilon\to0$. It was computed via the algorithm outlined in
Sect.~\ref{sec:pde}, with
$L=\frac1\alpha P \kappa \partial_x^2 - \textrm{Id}$ and
$R(u) = \frac1\alpha P(u-u^3)$. Starting at the fixed point $A$ the
minimizer takes a very different path than the heteroclinic orbit. It
first moves the domain wall, at vanishing cost for $\alpha\to0$,
without nucleation. At the point $X$ the motion changes, tracking
closely the separatrix towards the unstable point $S$.  From this point
onward, $S \to B$, the transition path then follows the heteroclinic
orbit, which is the deterministic relaxation path. In this respect,
the SPDE model \eqref{eq:acch-model} resembles closely the
2-dimensional model \eqref{eq:acch-toy}.

\begin{figure}[tb]
  \includegraphics[width=200pt]{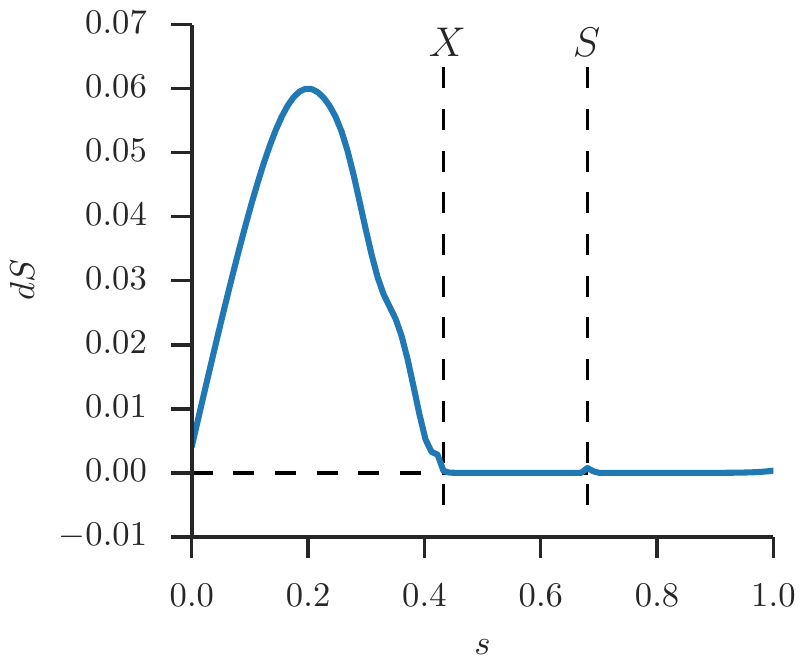}
  \caption{Action along the minimizer. Note that the action is
    non-zero climbing up the slow-manifold, but diminishes to zero
    already at $X$ when it approaches the separatrix, before it
    reaches~$S$. \label{fig:cornersb}}
\end{figure}

To further illustrate this resemblance, we choose to project the
minimizer and the heteroclinic orbit onto two coordinates,
\begin{enumerate}
\item its mean $\int \phi(x) \,dx$, which resembles the direction
  $\phi_1+\phi_2$ of the 2-dimensional model, and
\item its component in the direction of the initial (or final) state,
  $\int \phi(x) \phi_A(x)\,dx$, which corresponds to the direction
  $\phi_1-\phi_2$ of the 2-dimensional model.
\end{enumerate}
The transition path and the heteroclinic projected in these reduced
coordinates are depicted in Fig.~\ref{fig:projection}. Note that this
figure is not a schematic, but the actual projection of the
heteroclinic orbit and the minimizer of Fig.~\ref{fig:transition}
according to (i) and (ii) above. The separatrix is the straight line
$\int \phi(x)\phi_A(x)\,dx=0$. The movement of the minimizer (dark)
closely along the slow manifold (dashed), $A\to X$, and the
separatrix, $X\to S$, (which is also part of the slow manifold) into
$S$ highlights its difference with the heteroclinic orbit (light). The
configurations at the points $A, B, S$ and $X$ are depicted in
Fig.~\ref{fig:corners}, while Fig.~\ref{fig:cornersb} shows the action
density $dS$ along the transition path. Note that this quantity
becomes close to zero already at $X$, because the minimizer follows
closely the separatrix from $X$ to $S$, and this motion is therefore
quasi-deterministic.

The numerical parameters we used in these computations are
$h=10^{-1}$, $N_s = 100$, $N_x=2^6$, where $N_x$ denotes the number of
spatial discretization points.

\subsection{Burgers-Huxley model}
\index{Burgers-Huxley model}

As a second example involving an SPDE, we consider 
\begin{equation} \label{eq:burgers_general}
  u_t + \alpha uu_x - \kappa u_{xx} = f(u,x,t) + \sqrt{\epsilon} \eta(x,t)\,.
\end{equation}
where $\alpha>0$ and $\kappa>0$ are parameters, and we impose periodic
boundary condition on $x\in [0,1]$.  Without the term $f(u,x,t)$, this
is the stochastic Burgers equation which arises in a variety of
fields, in particular in the context of compressible gas dynamics,
traffic flow, and fluid dynamics. With the reaction term $f(u,x,t)$
added this equation is referred to as the (stochastic) Burgers-Huxley
equation \cite{wang-zhu-lu:1990} , which has been used e.g. to
describe the dynamics of neurons. The addition of a reaction term
makes it possible to obtain multiple stable fixed points. As a
particular case, we will consider (\ref{eq:burgers_general}) with
\begin{equation}
  \label{eq:burgers_react}
  f(u,x,t) = - u(1-u)(1+u)
\end{equation}
so that $u_+ = 1$ and $u_- = -1$ are the two stable fixed points of the
deterministic dynamics. We are interested in the mechanism of the
noise-induced transitions between these points.
\begin{figure}[tb]
  \begin{center}
    \includegraphics[width=0.53\textwidth]{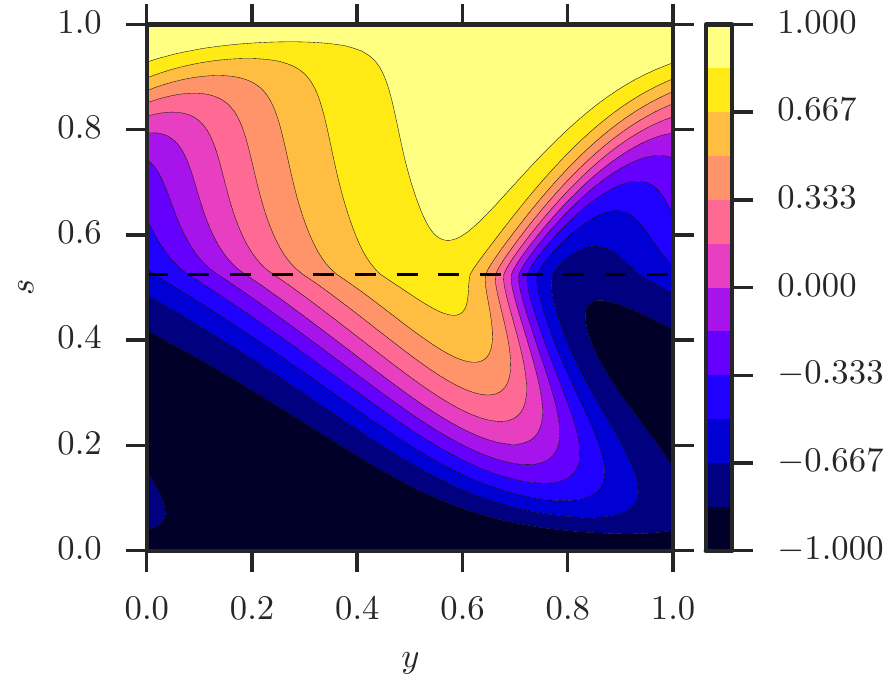}\hfill
    \includegraphics[width=0.45\textwidth]{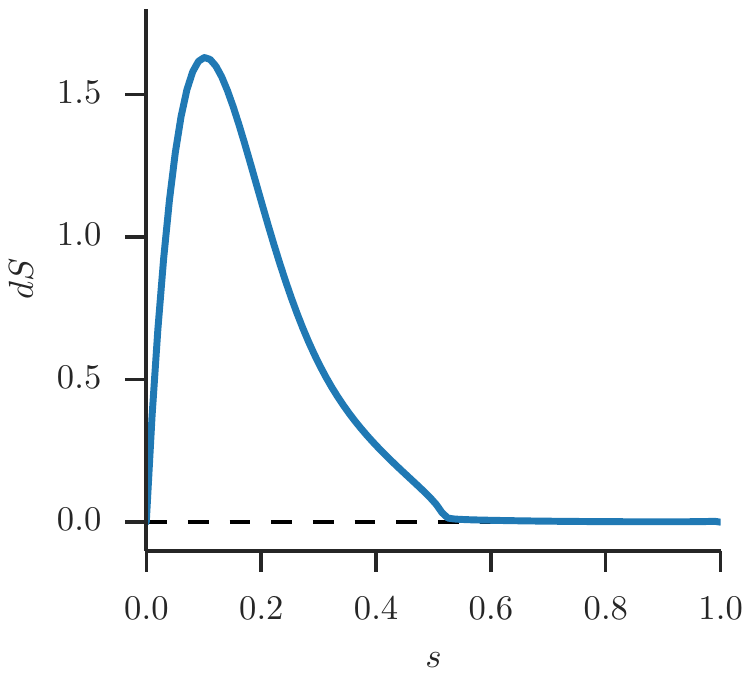}
  \end{center}
  \caption{Burgers-Huxley equation: Minimizer switching from $u_- =
    -1$ to $u_+ = 1$. Left: $u$-field. The saddle-point is marked with
    a dashed line. There is a noticeable kink in the dynamics
    switching from uphill ($s<s_{\text{saddle}}$) to downhill
    ($s>s_{\text{saddle}}$) dynamics. Right: Action density along the
    minimizer. \label{fig:BH-trajectory}}
\end{figure}
%
%

When $\alpha=0$, the system is in detailed balance and therefore the
forward and backward reaction follow the same path. The potential
associated with the reaction term~\eqref{eq:burgers_react} is
symmetric under $u\to-u$, and both states are equally probable. In
contrast, when $\alpha\ne0$ it is not obvious \textit{a priori}
whether $u_+$ and $u_-$ are equally probable, since the non-linearity
breaks the spatial symmetry, leading to a steepening of negative
gradients into shocks while flattening positive gradients. A
computation of the minimizer of the geometric action in both
directions, for $\kappa=0.01$ and $\alpha=\tfrac14$ reveals that
indeed forward and backward reactions are equally probable, even though
the transition paths do not coincide with the heteroclinic orbits. The
transition from $u_-$ to $u_+$ is depicted in Fig.
\ref{fig:BH-trajectory} (left). An intuitive explanation for the equal
probability of $u_+$ and $u_-$ is given by the fact that the backward
reaction pathways is identical to the forward path under the
transformation $u\to-u$, $x\to-x$. The action along this minimizer is
depicted in figure \ref{fig:BH-trajectory} (right). The minimizer is
computed via the algorithm lined out in Sect. \ref{sec:pde}, with
$L=-\kappa \partial_x^2$ and $R(u) = \alpha uu_x + u(1-u)(1+u)$.

The numerical parameters were chosen as $N_s = 100$, $N_x=2^8$,
$h=5\cdot 10^{-3}$.

\subsection{Noise-induced transitions between climate regimes\index{climate regimes}}
\label{ssec:meta-stable-climate}

Many climate systems exhibit metastability. Examples include the
Kuroshio oceanic current\index{Kuroshio oceanic current} off the coast
of Japan, which can be in either a small or a large meander state and
rarely switches between the two \cite{chao:1984,
  schmeits-dijkstra:2001}, or the atmospheric mid-latitude circulation
over the North-Atlantic, which makes rare transitions between a
strongly zonal and a weakly zonal (``blocked'') flow, characterized as
``Grosswetterlagen''\index{Grosswetterlagen} in \cite{baur:1951}. In
these and similar examples, the climate system stays trapped in the
vicinity of the stable regimes most of the time. Random noise,
originating either from physical stresses or from unresolved modes in
truncated models, induces rare regime transitions, which can be
captured by large deviation minimizers. The transition trajectory and
their corresponding action allow to make statements about not only the
relative probability of the different regimes and the transition
rates, but also the exact transition pathway taken to switch between
regimes.

We want to illustrate the feasibility of our numerical scheme for this
particular field of application by investigating metastability in two
simple climate models: A three-dimensional model for Grosswetterlagen
proposed by Egger \cite{egger:1981} and the six-dimensional
Charney-DeVore model \cite{charney-devore:1979}. Due to their highly
truncated nature, both models have very limited predictive power, but
exemplify the phenomenon of metastability in climate patterns or
regimes.

\subsubsection{Metastable climate regimes in Egger's model}
\index{Egger's model}

\begin{figure}[tb]
  \begin{center}
    \includegraphics[width=0.48\textwidth]{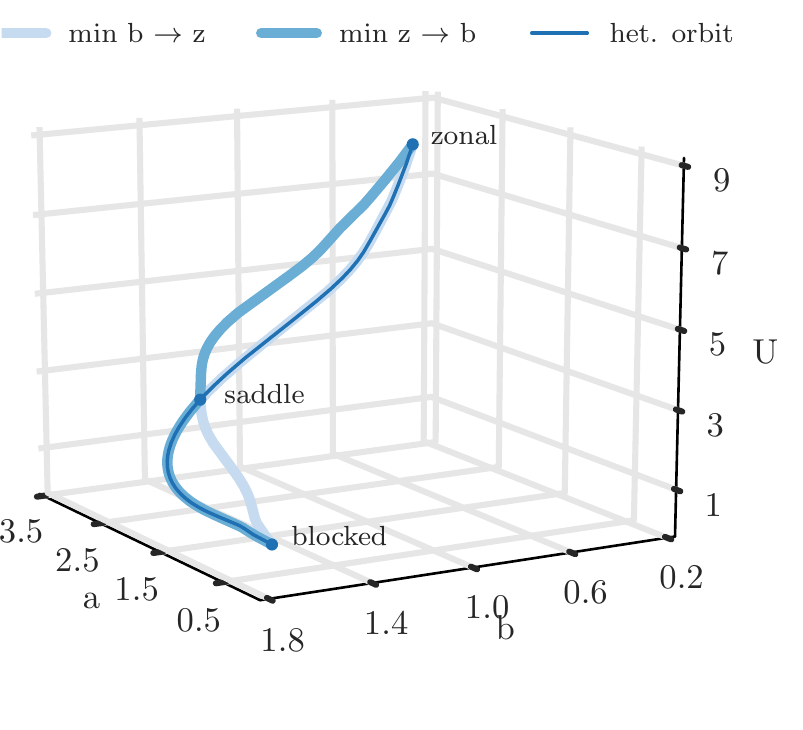}\hfill
    \includegraphics[width=0.48\textwidth]{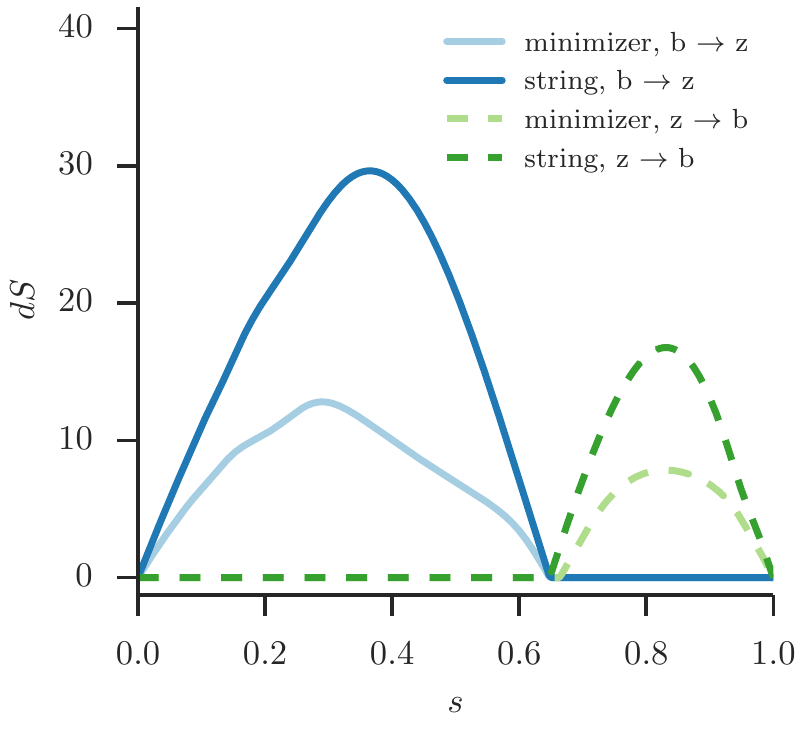}
  \end{center}
  \caption{Egger's model with $H = 12, \beta = 1.25, \gamma = 2, k = 2,
    U_0 = 10.5$ Left: Minimizers and deterministic relaxation
    paths. Right: Comparison of the action density.\label{fig:eggers}}
\end{figure}

Egger \cite{egger:1981} introduces the following SDE system as a
crude model to describe weather regimes in central Europe:
\begin{equation}
  \label{eq:eggers}
  \begin{cases}
    da = kb(U-\beta/k^2)\,dt - \gamma a\,dt + \sqrt{\epsilon} dW_a\,,\\
    db = -ka(U-\beta/k^2)\,dt + UH/k\,dt - \gamma b\,dt 
    + \sqrt{\epsilon} dW_b\,,\\
    dU = -bHk/2\,dt-\gamma(U-U_0)\,dt + \sqrt{\epsilon} dW_U\,.
  \end{cases}
\end{equation}
When $\epsilon$ is small, these equation exhibit metastability between
a ``blocked state'' and a ``zonal state'', shown in
Fig.~\ref{fig:eggers}. We use our gMAM algorithm to compute the
transition paths between these states. The system~\eqref{eq:eggers}
falls into the category discussed in Sect.~\ref{ssec:additive}
(i.e. additive noise) and can be solved with the simplest variant of
the algorithm. For $H = 12, \beta = 1.25, \gamma = 2, k = 2$ and $U_0
= 10.5$, the fixed points are approximately
$(a,b,U)=(0.465,1.65,0.593)$ for the blocked, $(3.07, 0.392, 8.15)$
for the zonal and $(2.80, 1.35, 2.38)$ for the unstable fixed point
(saddle). The minimizers of the action are show in
Fig.~\ref{fig:eggers}(left) where they are compared to the
heteroclinic orbits that connects the unstable critical points to the
stable ones. The action density along the transition trajectories and
the heteroclinic orbits is depicted in Fig.~\ref{fig:eggers} (right).

The numerical parameters we used in these computations are $N_s =
2^{8}$, $h=10^{-3}$.

\subsubsection{Metastable climate regimes in the Charney-DeVore model}
\index{Charney-DeVore model}

\begin{figure}[tb]
  \begin{center}
    \includegraphics[width=0.9\textwidth]{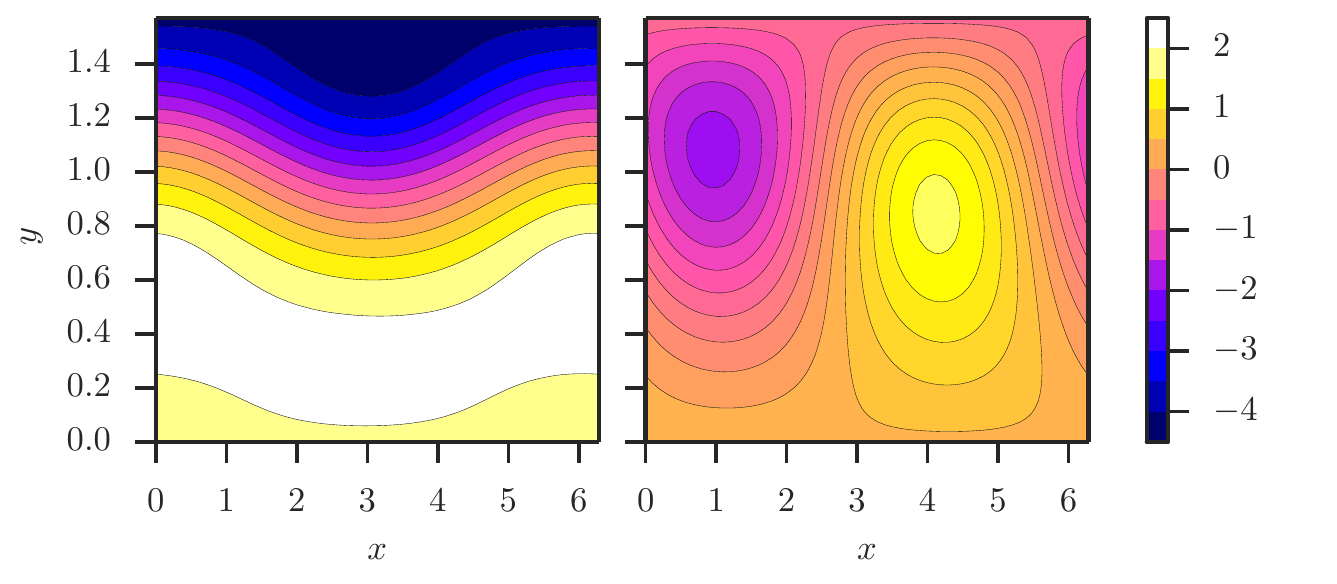}
  \end{center}
  \caption{Contours of the stream-function $\psi(x,y)$ of the two
    meta-stable configurations of the 6-dimensional CDV model. Left:
    Zonal state; Right: Blocked state.\label{fig:cdv-fixed-points}}
\end{figure}

\begin{figure}[tb]
  \begin{center}
    \includegraphics[width=0.96\textwidth]{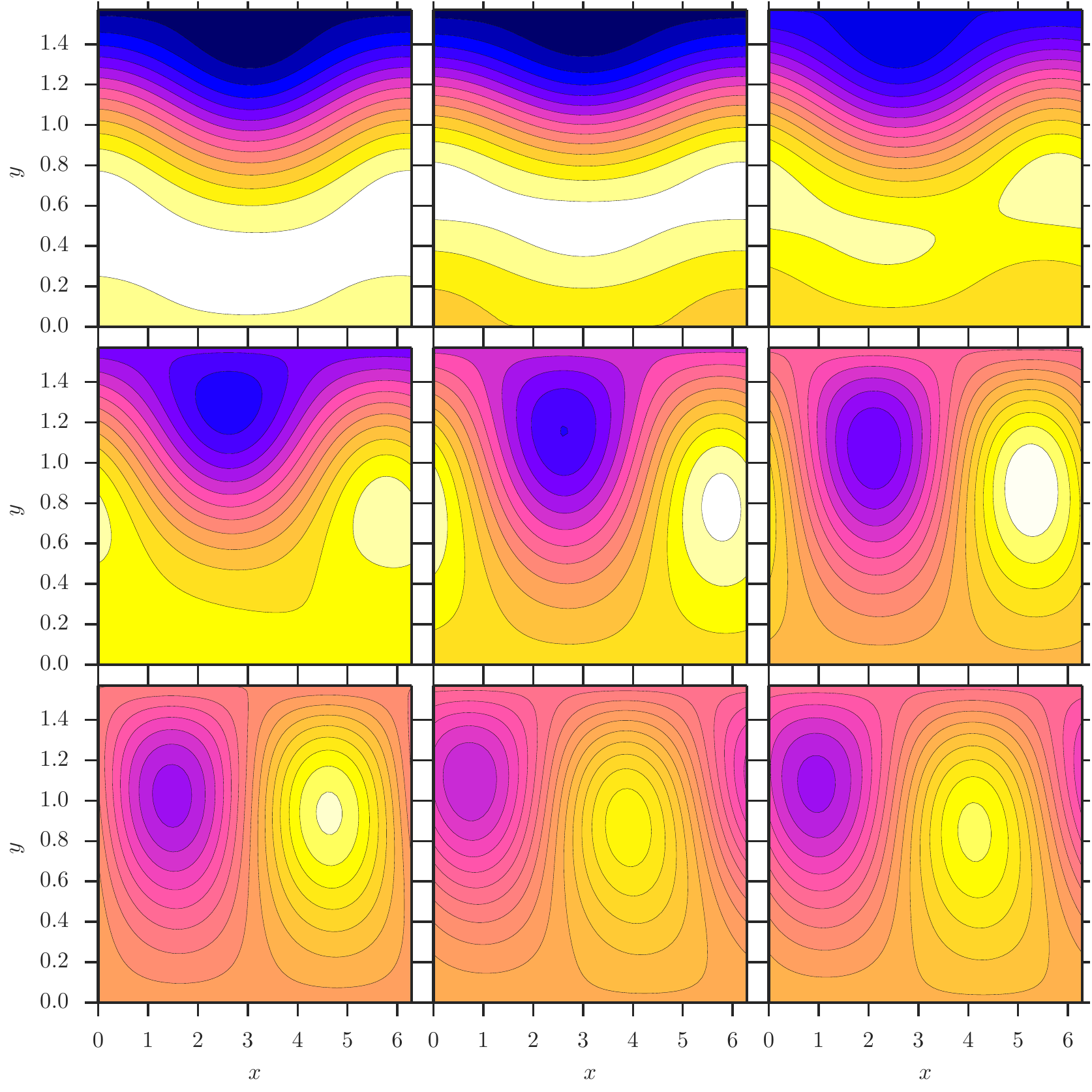}
  \end{center}
  \caption{Contours of the stream-function $\psi(x,y)$ along the
    transition trajectory from the zonal to the blocked meta-stable
    configuration for the CDV model. The arclength parameter increases
    in lexicographic order, with the top left plot being the initial
    state and the bottom right plot being the final state. The saddle
    point configuration is depicted in the center. The colormap is
    identical to figure
    \ref{fig:cdv-fixed-points}.\label{fig:cdv-zonal-to-blocked}}
\end{figure}

\begin{figure}[tb]
  \begin{center}
    \includegraphics[width=0.96\textwidth]{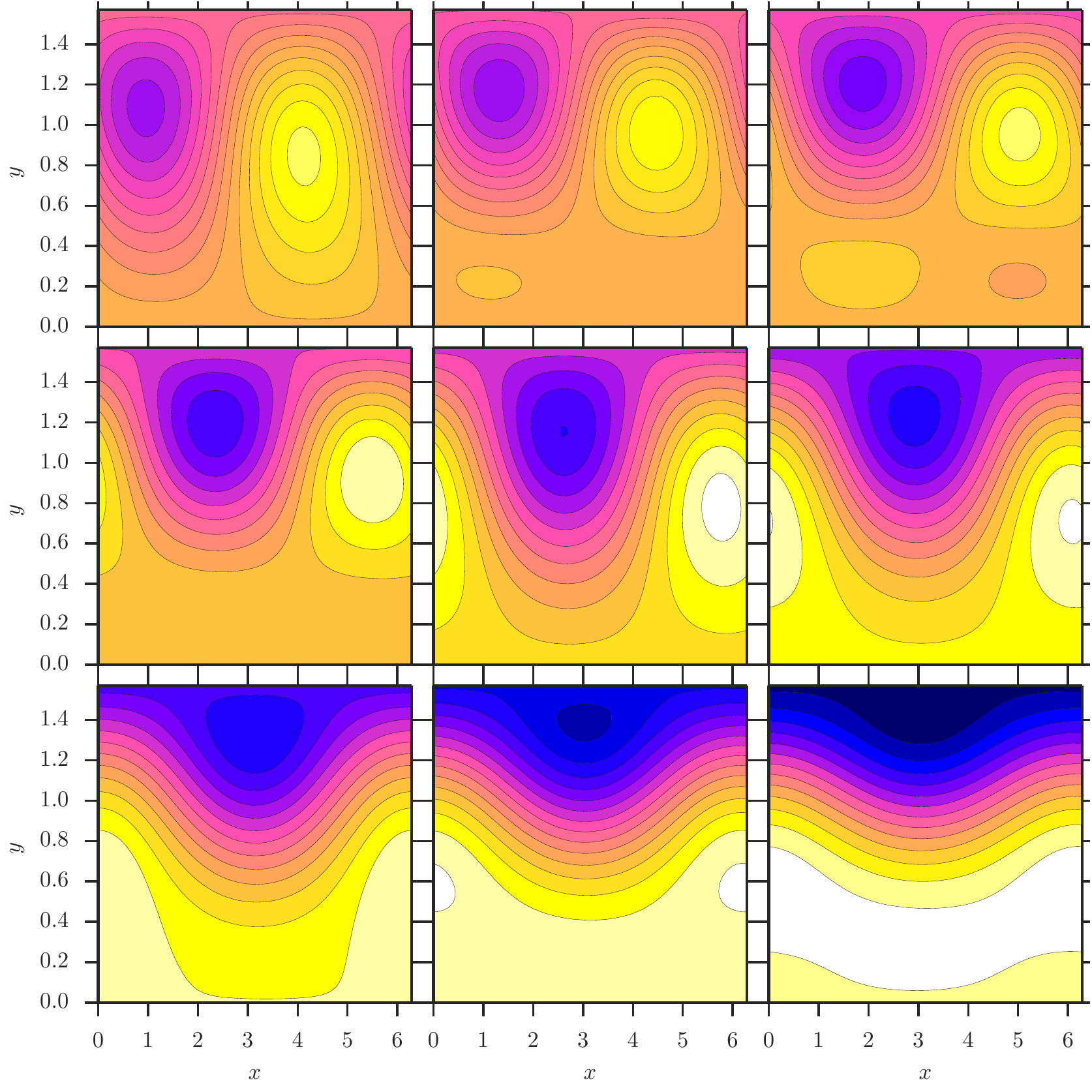}
  \end{center}
  \caption{Contours of the stream-function $\psi(x,y)$ along the
    transition trajectory from the blocked to the zonal meta-stable
    configuration for the CDV model. The arclength parameter increases
    in lexicographic order, with the top left plot being the initial
    state and the bottom right plot being the final state. The saddle
    point configuration is depicted in the center. The colormap is
    identical to figure
    \ref{fig:cdv-fixed-points}.\label{fig:cdv-blocked-to-zonal}}
\end{figure}

\begin{figure}[tb]
  \includegraphics[width=160pt]{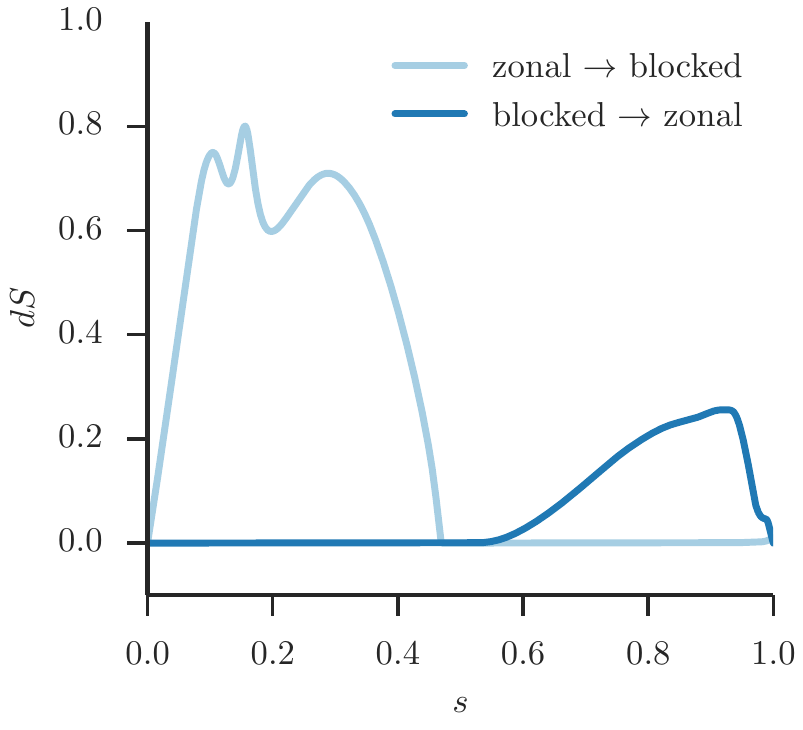}
  \caption{Action density $dS$ along the transition pathways from
    zonal to blocked (forward) and from blocked to zonal
    (backward). In both directions, after passing the saddle point,
    the action becomes zero since the motion is
    deterministic. \label{fig:cdv-action}}
\end{figure}

Egger's model retains no nonlinear interaction between different fluid
modes, which is believed to be insufficient to explain the transitions
between zonal and blocked states. A more sophisticated model,
truncating the barotropic vorticity equation\index{barotropic
  vorticity equation} (BVE) with full nonlinear terms, was introduced
by Charney and DeVore \cite{charney-devore:1979}. Their starting point
is the two-dimensional BVE on the $\beta$-plane,
\begin{equation}
  \label{eq:vorticity}
  \frac\partial{\partial t} \omega = u\cdot\nabla\omega - C(\omega-\omega^*)\,.
\end{equation}
Here $\omega = \zeta + \beta y + \gamma h$ is the total vorticity,
where $\gamma h$ is the topography in the $\beta$-plane, with $\beta =
2\Omega \cos(\theta)/R$ for planetary angular velocity $\Omega$,
radius $R$ and latitude $\theta$, and $\zeta = \Delta \psi$ is the
relative vorticity for the stream-function\index{stream-function}
$\psi$. The term $-C(\omega-\omega^*)$ accounts for Ekman
damping\index{Ekman damping} with coefficient $C>0$.

Charney-DeVore considered the vorticity equation~\eqref{eq:vorticity}
in the box $[0,2\pi]\times[0,\pi b]$ with periodic boundary conditions
in $x$-direction and no-slip boundary conditions in
$y$-direction. They then projected this equation over 6 Fourier modes
in total, using the following representation for the stream-function
$\psi(x,y,t)$:
\begin{equation}
  \psi(x,y,t) = \sum_{n,m} \psi_{nm}(t)\phi_{nm}(x,y)\,,
\end{equation} 
where the sums run on $n\in\{-1,0,1\}$ and $m\in\{1,2\}$ and
\begin{align}
  \phi_{0m}(y) &= \sqrt{2} \cos(my/b),& \phi_{nm}(x,y) &= \sqrt{2} e^{inx} \sin(my/b)\,.
\end{align}
Letting $x_i$, $i\in\{1,\dots,6\}$ be defined as
\begin{equation}
  \begin{split}
    x_1 = \frac1b \psi_{01}, \qquad x_2 = \frac{1}{\sqrt{2}b} \left( \psi_{11} + \psi_{-11}\right), \qquad x_3 = \frac{i}{\sqrt{2}b} \left(\psi_{11} - \psi_{-11}\right),\\
    x_4 = \frac1b \psi_{02}, \qquad x_5 = \frac{1}{\sqrt{2}b} \left( \psi_{12} + \psi_{-12}\right), \qquad x_6 = \frac{i}{\sqrt{2}b} \left(\psi_{12} - \psi_{-12}\right)\,,
  \end{split}
\end{equation}
taking the following form for the topography
\begin{equation}
  h(x,y) = \cos(x)\sin(y/b)\,,
\end{equation}
and choosing $\omega^*$ such that only two parameters $x_1^*$ and
$x_4^*$ are free and the other are set zero, they arrived at the
following six-dimensional model
\begin{equation}
  \begin{split}
    d x_1 &= \left(\tilde \gamma_1 x_3 - C(x_1-x_1^*)\right) dt + \sqrt{2\epsilon} dW_1\,,\\
    d x_2 &= \left(-(\alpha_1x_1 - \beta_1)x_3 - Cx_2 -\delta_1 x_4x_6\right) dt + \sqrt{2\epsilon} dW_2\,,\\
    d x_3 &= \left((\alpha_1x_1-\beta_1)x_2 - \gamma_1 x_1 - Cx_3 + \delta_1x_4x_5\right) dt + \sqrt{2\epsilon} dW_3\,,\\
    d x_4 &= \left(\tilde\gamma_2 x_6-C(x_4-x_4^*)+\eta(x_2x_6-x_3x_5)\right) dt + \sqrt{2\epsilon} dW_4\,,\\
    d x_5 &= \left(-(\alpha_2x_1-\beta_2)x_6-Cx_5-\delta_2x_3x_4\right) dt + \sqrt{2\epsilon} dW_5\,,\\
    d x_6 &= \left((\alpha_2x_1-\beta_2)x_5 - \gamma_2x_4 - Cx_6 + \delta_2x_2x_4\right) dt + \sqrt{2\epsilon} dW_6\,,
  \end{split}
\end{equation}
where, for $m\in\{1,2\}$,
\begin{equation}
  \begin{split}
    \alpha_m &= \frac{8\sqrt{2}}{\pi} \frac{m^2}{4m^2-1}\frac{b^2+m^2-1}{b^2+m^2}\,,\\
    \beta_m &= \frac{\beta b^2}{b^2+m^2}\,,\\
    \gamma_m &= \gamma \frac{\sqrt{2}b}{\pi}\frac{4m^3}{(4m^2-1)(b^2+m^2)}\,,\\
    \tilde\gamma_m &= \gamma \frac{\sqrt{2}b}{\pi} \frac{4m}{4m^2-1}\,,\\
    \delta_m &= \frac{64\sqrt{2}}{15\pi} \frac{b^2-m^2+1}{b^2+m^2}\,,\\
    \eta &= \frac{16\sqrt{2}}{5\pi}\,.
  \end{split}
\end{equation}
The original Charney-DeVore equation did not contain random forcing
terms: here we added to each equations an independent white noise
$dW_i$ with amplitude $\sqrt{2\epsilon}$.

Choosing $b=\tfrac12$, $C = \tfrac1{10}$, $\beta = \tfrac54$, $\gamma
= 1$, $x_1^* = \tfrac92$, and $ x_4^* = -\tfrac95$, the 6-dimensional
stochastic model above possesses two metastable states, shown in
Fig.~\ref{fig:cdv-fixed-points}: a zonal state (left) and a blocked
state (right). The transition paths from zonal to blocked and from
blocked to zonal are different. They are shown in
Fig.~\ref{fig:cdv-zonal-to-blocked} and
Fig.~\ref{fig:cdv-blocked-to-zonal}, respectively, and they were both
calculated by minimizing the geometric action using our simplified
gMAM algorithm. The actions along both paths are depicted in
Fig.~\ref{fig:cdv-action}.

The numerical parameters in these computations were $N_s = 100$,
$h=10^{-3}$.

\subsection{Generalized voter/Ising model}
\index{voter model}
\index{Ising model}
\label{ssec:voter-model}

\begin{figure}[tb]
  \begin{center}
    \includegraphics[width=0.48\textwidth]{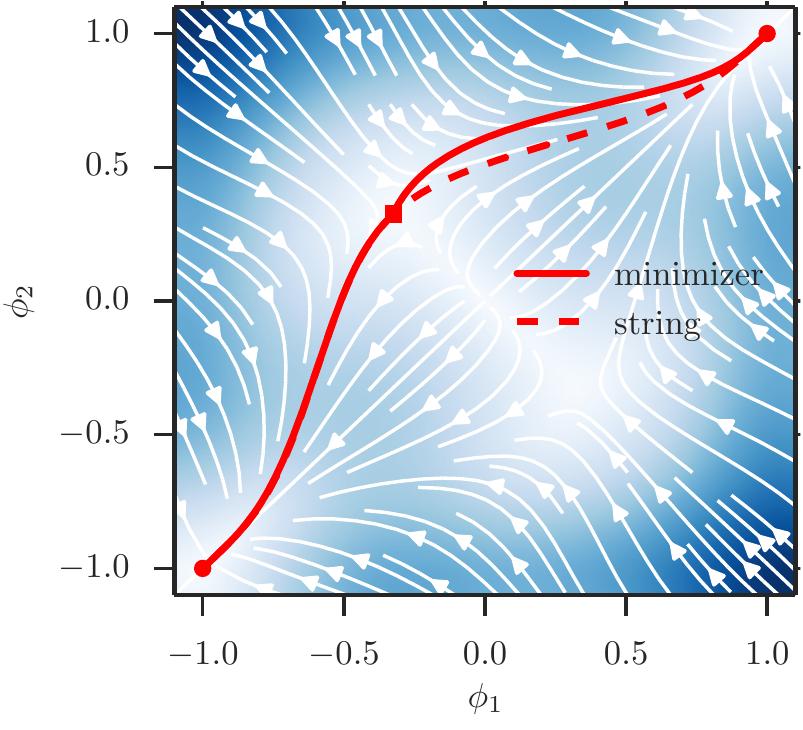}\hfill
    \includegraphics[width=0.48\textwidth]{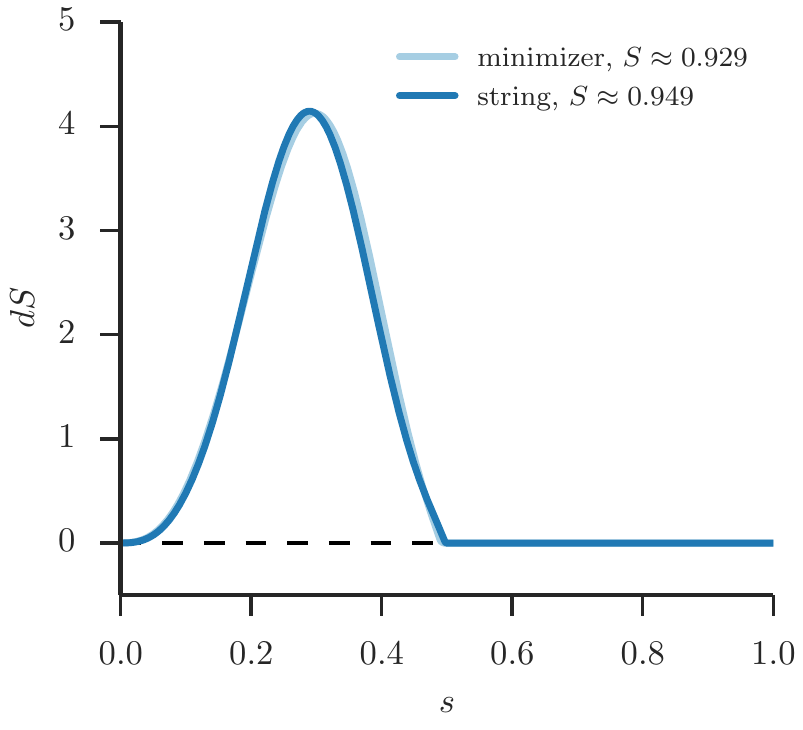}
  \end{center}
  \caption{Generalized voter/Ising model. Left: The arrows denote the
    direction of the deterministic flow, the shading its
    magnitude. The solid line depicts the minimizer, the dashed line
    the heteroclinic orbit. Markers are located at the fixed
    points (circle: stable; square: saddle). Right: Action density
    along the minimizers for the two trajectories, with normalized
    path parameter $s\in(0,1)$.\label{fig:generalized-voter}}
\end{figure}

To analyze phase transitions in out-of-equilibrium systems, a Langevin
equation was proposed in \cite{al_hammal-chate-dornic-etal:2005} that
models critical phenomena with two absorbing states\index{absorbing
  state}. This equation was constructed by requiring that it be
symmetric under the transformation $\phi\to-\phi$ and have two
absorbing states, arbitrarily chosen to be at $\pm1$. The presence of
these absorbing states makes the noise multiplicative, with a scaling
involving the square root of the distance to the absorbing boundaries,
as suggested by the voter model \cite{dickman-tretyakov:1995,
  dornic-chate-chave-etal:2001}.  In order to account for Ising-like
spontaneous symmetry breaking\index{spontaneous symmetry breaking},
the authors of \cite{al_hammal-chate-dornic-etal:2005} also added a
bi-stable ``potential''-term with $-V'(\phi) = (a\phi-b\phi^3)$ to the
equation, which finally lead them to:
\begin{equation}
  \label{eq:voter_pde}
  \phi_t = \left((1-\phi^2)(a\phi-b\phi^3) + D \phi_{xx} \right)\,dt +
  \sigma \sqrt{1-\phi^2} \eta(x,t)\,.
\end{equation}
In the absence of noise ($\epsilon=0$) and for $a>0$, the $\phi=0$
state is locally unstable, but $b>0$ ensures stable fixed points at
$\phi=\pm \sqrt{a/b}$. In the limit $a/b\to1$, these fixed points
approach the absorbing boundaries, and we are interested in the noise
induced transition between these states. 

We stress that making mathematical sense of \eqref{eq:voter_pde} is
non-trivial (see the discussion in Sect.~\ref{sec:pde}). In the present
application, we are going to consider a finite truncation of this
SPDE, where the question of spatial regularity
disappears. Specifically, we transform~\eqref{eq:voter_pde} into a
two-dimensional stochastic ODE model by discretizing the spatial
direction via the standard 3-point Laplace stencil, and taking only
$N_x=2$ discretization points. This yields the stochastic ODE system
\begin{equation}
  \label{eq:generalized-voter}
  \begin{cases}
    d\phi_1 &= \left((1-\phi_1^2)(a\phi_1-b\phi_1^3) + D (\phi_1 - \phi_2)\right)\,dt + \sigma
    \sqrt{1-\phi_1^2} \,dW_x\\ 
    d\phi_2 &= \left((1-\phi_2^2)(a\phi_2-b\phi_2^3) - D (\phi_1 - \phi_2) \right)\,dt + \sigma
    \sqrt{1-\phi_2^2} \,dW_y\,,    
  \end{cases}
\end{equation}
where the constant $D$ couples the two degrees of freedom. This SDE
poses an interesting test-case for our numerical scheme, since not
only the noise is multiplicative\index{multiplicative noise}, but also
the computational domain must be restricted. The square defined by
$1=\max(|\phi_1|,|\phi_2|)$ marks the region in which the noise is
defined (real), and the noise decreases towards zero as it approaches
this absorbing barrier. Analog to the discussion in
\cite{al_hammal-chate-dornic-etal:2005}, the choice of the parameters
$(a,b)$ determines the dynamics, in particular if $a>0,b>0$ the model
exhibits bi-stability: There is an unstable fixed point at
$\phi=(0,0)$ and stable fixed points at $\phi = \pm
(\sqrt{a/b},\sqrt{a/b})$. As long as $a<b$, these fixed points are
inside the allowed region. For $a/b \to 1$ the two stable fixed points
approach the absorbing boundary. Here, we take $b=1, a=1-10^{-4},
D=0.4$, so that $\sqrt{a/b} \approx 0.99995$ is located close to the
barrier at $1$. The minimizer and corresponding action are shown in
Fig. \ref{fig:generalized-voter}.

The numerical parameters were chosen as $N_s = 2^{8}$, $h=10^{-3}$.

\subsection{Bi-stable reaction-diffusion model}
\index{reaction-diffusion equation}

In the context of chemical reactions and birth-death
processes\index{birth-death process}, one considers networks of
several reactants in a container of volume $V$ which is considered
well-stirred. As an example case, we consider the bi-stable chemical
reaction network\index{chemical reaction network}
\begin{equation*}
  A \underset{k_1}{\stackrel{k_0}{\rightleftharpoons}} X,\qquad 2X+B \underset{k_3}{\stackrel{k_2}{\rightleftharpoons}} 3X
\end{equation*}
with rates $k_i>0$, and where the concentrations of $A$ and $B$ are
held constant. This system was introduced in \cite{schloegl:1972} as a
prototypical model for a bi-stable reaction network. Its dynamics can
be modeled as a Markov jump process (MJP)\index{Markov jump process}
with generator
\begin{equation}
  \label{eq:prob22}
  (L^Rf)(n) = A_+(n) \left(f(n+1) - f(n)\right) + A_-(n)\left(f(n-1) - f(n)\right)
\end{equation}
with the propensity functions
\begin{equation}
  \begin{cases}
    A_+(n) &= k_0 V + (k_2/V) n(n-1)\\
    A_-(n) &= k_1 n + (k_3/V^2) n(n-1)(n-2)\,.
  \end{cases}
\end{equation}

The model above satisfies a large deviation principle in the following
scaling limit: Denote by $c=n/V$ the concentration of $X$, and
normalize it by a typical concentration, $\rho=c/c_0$. Now, in the
limit of a large number of particles per cell $\Omega=c_0 V$ and
simultaneously rescaling time by $1/\Omega$, we obtain
\begin{equation}
  \label{eq:prob22scaled}
  (L_\epsilon^Rf)(\rho) = \frac1\epsilon \Big(a_+(\rho)
  \left(f(\rho+\epsilon)
    -f(\rho)\right) +
  a_-(\rho)\left(f(\rho-\epsilon)-f(\rho)\right)\Big)\,,
\end{equation}
where $\epsilon=1/\Omega$ is a small parameter. Here, we defined
$k_i=\lambda_i(c_0)^{1-i}$, and
\begin{equation}
  \begin{cases}
    a_+(\rho) &= \lambda_0 + \lambda_2 \rho^2\\
    a_-(\rho) &= \lambda_1 \rho + \lambda_3 \rho^3\,.
  \end{cases}
\end{equation}
The large deviation principle for~\eqref{eq:prob22scaled} can be
formally obtained via WKB analysis\index{WKB approximation}, that is,
by setting $f(\rho) = e^{\epsilon^{-1}G(\rho)}$
in~\eqref{eq:prob22scaled} and expanding in $\epsilon$
\cite{doering-sargsyan-sander-etal:2007}. To leading order in
$\epsilon$, this gives an Hamilton-Jacobi
operator\index{Hamilton-Jacobi operator} associated with an
Hamiltonian that is also the one rigorously derived in
LDT~\cite{shwartz-weiss:1995}. It reads
\begin{equation}
  \label{eq:H-reactive}
  H(\rho,\vartheta) = a_+(\rho)(e^\vartheta-1) + a_-(\rho)
  (e^{-\vartheta}-1)\,.
\end{equation}
This is an example of a system whose Hamiltonian is not quadratic in
the conjugate momentum $\vartheta$\index{non-quadratic
  Hamiltonian}. Therefore the computation of $\vartheta_*$ by
\eqref{eq:theta-star} can not be performed explicitly in general. For
parameters $\lambda_0 = 0.8, \lambda_1 = 2.9, \lambda_2 = 3.1,
\lambda_3 = 1$, the system has two stable fixed points $\rho_{\pm}$
and a saddle $\rho_s$ at $\rho_+=\tfrac85, \rho_-=\tfrac12, \rho_s=1$.

Since transitions in 1D are fairly trivial, we want to consider the
case of $N$ neighboring reaction compartments, each well-stirred, but
with random jumps possible between neighboring compartments. This
situation was analyzed in \cite{tanase-nicola-lubensky:2012} via
direct sampling, but we are interested in the computation of the
transition trajectory. Denote by $\rho_i$ the concentration in the
$i$-th compartment and refer to the vector $\pmb{\rho}$ as the
complete state, $\pmb{\rho} = \sum_{i=0}^N \rho_i \hat e_i$. In this
case, we obtain a diffusive part of the generator, $L^D$, coupling
neighboring compartments. For a diffusivity $D$, it is
\begin{equation}
  (L^D f)(\pmb{\rho}) = \frac{D}{\epsilon} \sum_{i=1}^{N} \rho_i \left(f(\pmb{\rho} - \epsilon \hat e_i + \epsilon \hat e_{i-1}) + f(\pmb{\rho} - \epsilon \hat e_i + \epsilon \hat e_{i+1})- 2f(\pmb{\rho})\right)\,.
\end{equation}
The process associated with this generator also admits a large
deviation principle with Hamiltonian
\begin{equation}
  \label{eq:Hd1}
  H^D(\pmb{\rho},\pmb{\vartheta}) = D\sum_{i=1}^N \rho_i
  \left(e^{\vartheta_{i-1}-\vartheta_i} +
  e^{\vartheta_{i+1}-\vartheta_i} - 2\right)\,.
\end{equation}
Therefore, the full Hamiltonian becomes
$H(\pmb{\rho},\pmb{\vartheta}) = H^D(\pmb{\rho},\pmb{\vartheta}) +
\sum_{i=1}^N H^R(\rho_i,\vartheta_i)$,
where $H^R(\rho_i,\vartheta_i)$ is the reactive Hamiltonian in
\eqref{eq:H-reactive}, which is summed up over all the compartments.

\begin{figure}[tb]
  \begin{center}
    \includegraphics[width=265pt]{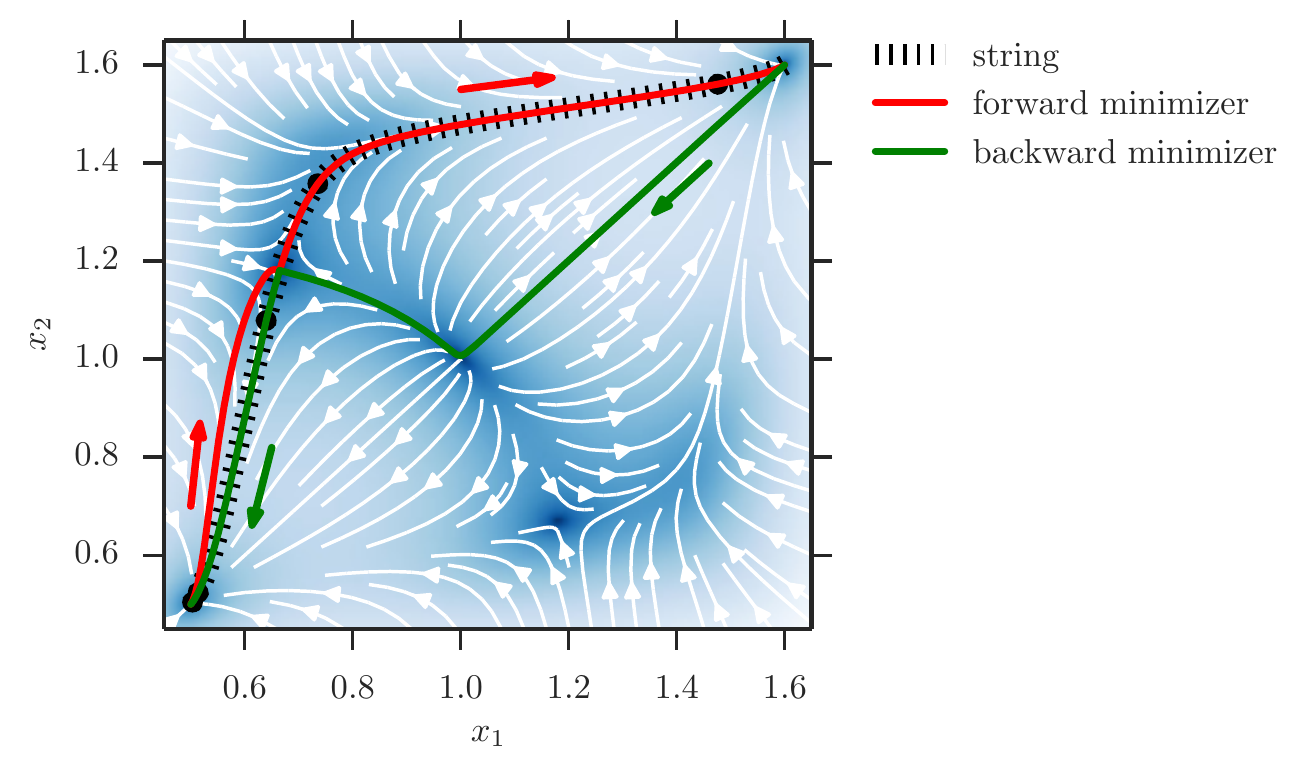}
  \end{center}
  \caption{Bi-Stable reaction-diffusion model with $N=2$ reaction
    cells. Show are the forward (red) and backward (green) transitions
    between the two stable fixed points, in comparison to the
    heteroclinic orbit (dashed). The flow-lines depict the
    deterministic dynamics, their magnitude is indicated by the
    background shading.\label{fig:reaction-diffusion}}
\end{figure}

\begin{figure}[t]
  \includegraphics[width=160pt]{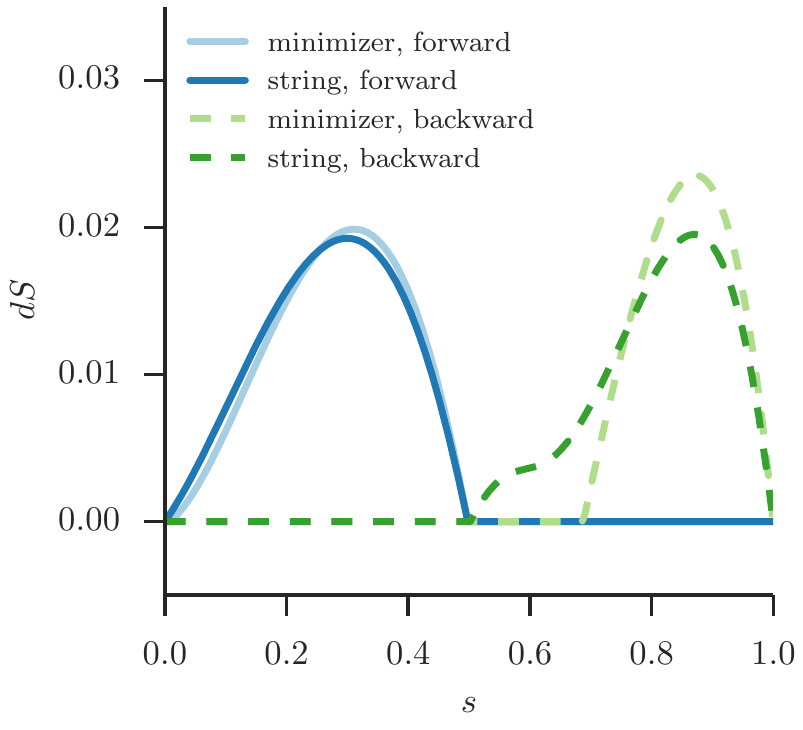}
  \caption{Action densities for the bi-stable reaction-diffusion
    model. Depicted are the actions corresponding to the forward
    (solid) and backward (dashed) minimizer (dark) and heteroclinic (light)
    orbit.\label{fig:reaction-diffusion-action}}
\end{figure}

We used our new gMAM algorithm to minimize the geometric action and
compute the transition paths between the stable fixed points for the
simplest non-trivial case of $N=2$ compartments. Shown in
Fig.~\ref{fig:reaction-diffusion} are the forward and backward
trajectories. Note that the backward transition
($(\rho_+,\rho_+)\to (\rho_-,\rho_-)$) takes a special form: It climbs
against the deterministic dynamics up to the maximum, then relaxes
along the separatrix down to the saddle. Additionally, we compare
these trajectories with the heteroclinic orbit obtained by the string
method. The action along these trajectories is depicted in
Fig.~\ref{fig:reaction-diffusion-action}. Note how for the backward
minimizer the action is zero already before it hits the saddle, as the
movement from the maximum to the saddle happens deterministically.

The numerical parameters were chosen as $N_s = 2^{9}$, $h=10$.

\subsection{Slow-fast systems}
\index{slow-fast system}

In contrast to a large deviation principle arising in the limit of
small noise or large number of particles, a different class of
Hamiltonians arises for systems with a slow variable $X$ evolving on a
timescale $O(1)$ and a fast variable $Y$ on a time scale $O(\alpha)$:
\begin{subequations}
  \begin{align}
    \dot X &= f(X,Y)\\
    dY &= \frac1\alpha b(X,Y)dt + 
         \frac1{\sqrt{\alpha}} \sigma(X,Y)dW\label{eq:fastslow-Y}\,.
  \end{align}
\end{subequations}
Examples of systems with large timescale separation\index{timescale
  separation} $\alpha\ll1$ are ubiquitous in nature, and usually one
is interested mostly in the long-time behavior of the slow
variables. In particular, we are concerned with situations where the
slow dynamics exhibits metastability. We want to use our algorithm to
compute transition pathways in this setup for the limit of infinite
time scale separation.

In the limit as $\alpha\to0$, the fast variables reach statistical
equilibrium before any motion of the slow variables, and these slow
variables only experience the average effect of the slow ones. This
behavior can be captured by the following deterministic limiting
equation which is akin to a law of large numbers\index{law of large
  numbers} (LLN) in the present context and reads
\begin{equation}
  \label{eq:fastslow-lln}
  \Dot{\bar X} = F(\bar X) \quad \text{where} \quad F(x) = 
  \lim_{T\to\infty} \frac1T \int_0^T f(x, Y_{x}(\tau))\,d\tau\,.
\end{equation}
Here $Y_x(t)$ is the solution of \eqref{eq:fastslow-Y} for $X(t)=x$
fixed \cite{freidlin:1978, bakhtin:2003, kifer:2009,
  bouchet-grafke-tangarife-etal:2016}. For small but finite $\alpha$,
the slow variables also experience fluctuations through the fast
variables. In particular, the statistics of
$\xi=(X-\bar X)/\sqrt{\alpha}$ on $O(1)$ time scales can be described
by a central limit theorem (CLT) as small Gaussian noise on top of the
slow mean $\bar X$. The CLT scaling, however, is inappropriate to
describe the fluctuations of the slow variables that are induced by
the effect of the fast variables on longer time scales and may, for
example, lead to transitions between stable fixed points of the
limiting equation in~\eqref{eq:fastslow-lln}. In particular, the naive
procedure of constructing an SDE out of the LLN and CLT to then
compute its LDT fails. Instead, the transitions in the limit of
$\alpha\to0$ are captured by an LDP with the Hamiltonian
\begin{equation}
  H(x,\vartheta) = \lim_{T\to\infty} \frac1T \log \mathbb{E} \exp\left(\vartheta
    \int_0^T f(x,Y_x(t))\,dt\right)\,.
  \label{eq:slowfast-H}
\end{equation}
Except for the special case $f(x,y) = r(x) + s(y)y$ (linear dependence
on the fast variable), the Hamiltonian \eqref{eq:slowfast-H} is
non-quadratic in $\theta$. As a consequence no S(P)DE with Gaussian
noise exists for the slow variable which has an LDP to describe the
transitions correctly.

The implicit nature of the Hamiltonian~\eqref{eq:slowfast-H}, in
particular containing an expectation, complicates numerical procedures
to compute its associated minimizers. Yet, in the non-trivial case of
a quadratic dependence of the slow variable on the fast ones, for
example,
\begin{equation}
  \begin{cases}
    \dot X = Y^2 - \beta X\\
    \displaystyle dY = -\frac1{\alpha}\gamma(X)Y\,dt 
    + \frac{\sigma}{\sqrt{\alpha}} dW\,,
  \end{cases}
\end{equation}
one indeed does obtain an explicit formula for the Hamiltonian
\eqref{eq:slowfast-H} (as derived in
\cite{bouchet-grafke-tangarife-etal:2016})
\begin{equation}
  \label{eq:H-slow-fast}
  h(x,\vartheta) = -\beta x \vartheta + \tfrac12\left(\gamma(x) 
    - \sqrt{\gamma^2(x) - 2\sigma^2 \vartheta}\right)\,.
\end{equation}
This example is interesting for our purpose not only because the
Hamiltonian is non-quadratic, but furthermore because of the existence
of a forbidden region $\vartheta > \gamma^2/(2\sigma)$ where the
Hamiltonian is not defined.

\begin{figure}[tb]
  \begin{center}
    \includegraphics[width=0.48\textwidth]{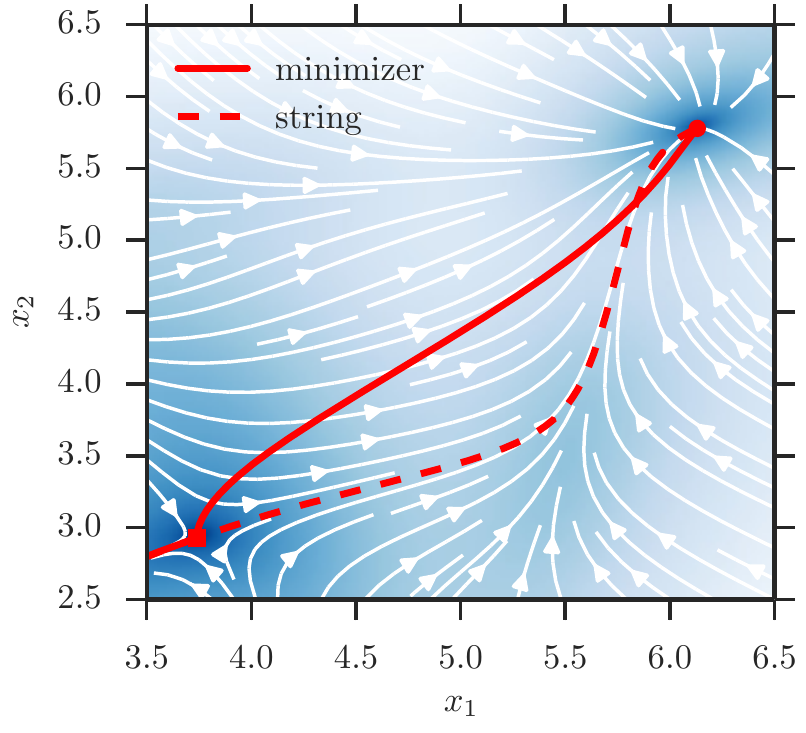}\hfill
    \includegraphics[width=0.48\textwidth]{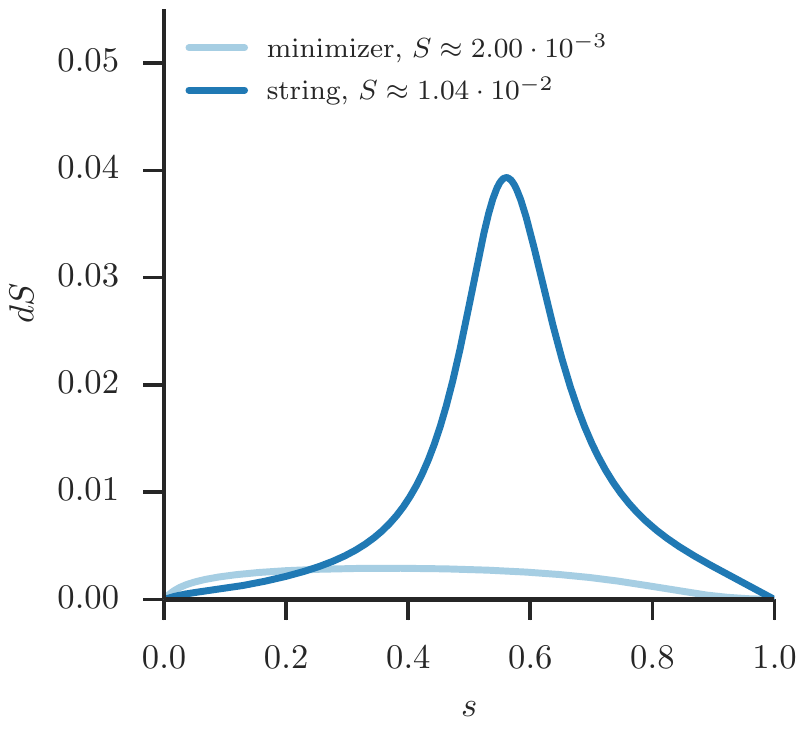}
  \end{center}
  \caption{Coupled slow-fast system ODE model for $D=1.0$. Left: The
    arrows denote the direction of the deterministic flow, the shading
    its magnitude. The solid line depicts the minimizer, the dashed
    line the relaxation paths from the saddle. Markers are located at
    the fixed points (circle: stable; square: saddle). Right: Action
    density along the minimizers for the two trajectories up to the
    saddle, with normalized path parameter
    $s\in(0,1)$.\label{fig:slowfast2}}
\end{figure}

Additionally increasing the number of degrees of freedom by combining
two independent multi-stable slow-fast systems and coupling them by a
spring with spring constant $D$, the full system reads
\begin{equation}
  \begin{cases}
    \dot X_1 = Y_1^2 - \beta_1 X_1 - D (X_1 - X_2)\\
    \dot X_2 = Y_2^2 - \beta_2 X_2 - D (X_2 - X_1)\\
   \displaystyle  dY_1 = -\frac1 {\alpha}\gamma(X_1)Y_1dt 
   + \frac{\sigma}{\sqrt{\alpha}} dW_1\\
    \displaystyle  dY_2 = -\frac1 {\alpha}\gamma(X_2)Y_2dt 
    + \frac{\sigma}{\sqrt{\alpha}} dW_2\,.
  \end{cases}
\end{equation}
The Hamiltonian for the LDT for this system is
\begin{equation}
  \label{eq:H-slowfast2}
  H(x_1,x_2,\vartheta_1,\vartheta_2) = h(x_1,\vartheta_1) 
  + h(x_2,\vartheta_2) + \langle -\nabla U(x_1,x_2), \vartheta \rangle\,,
\end{equation}
for $U(x,y) = \frac12D (x-y)^2$ and $h(x,\vartheta)$ defined as in
equation \eqref{eq:H-slow-fast}. The choice $\gamma(X)=(X-5)^2+1$
ensures two stable fixed points. The deterministic dynamics of this
system (i.e. the evolution of the averaged slow variables) are
depicted as white arrows in Fig. \ref{fig:slowfast2} (left). To
stress the important portion of the transition trajectory, the plot is
focused only on the initial state up to the saddle. Compared are the
minimizer and the heteroclinic orbits connecting the stable fixed
points to the saddle point. The corresponding actions are shown in
Fig. \ref{fig:slowfast2} (right). The specific choice of model
parameters for this computation is $\beta_1=0.6, \beta_2=0.3, D=1.0$
and $\sigma^2=10$.

The numerical parameters were chosen as $N_s = 2^{10}$, $h=10^{-2}$.

\section{Concluding Remarks}

We have discussed numerical schemes to compute minimizers of large
deviation action functionals, which are based on the geometric minimum
action method. The basis of these schemes is the minimization of a
geometric action on the space of arc-length parametrized curves, which
makes it possible to perform the double minimization over transition
time $T$ and action $S_T$ that is required to compute the LDT
quasipotential. In particular, transitions between metastable fixed
points of a system, which generally involve $T\to\infty$ and which are
not tractable with non-geometric minimum action methods can be
naturally analyzed in this setup.

A simplified gMAM algorithm was proposed here which is based on a
particular formulation of the geometric action leading to a mixed
optimization problem. This new formulation of the gMAM algorithm is
easier to implement than the original method: In its simplest form,
only first order derivatives of the Hamiltonian
$H(\varphi, \vartheta)$ are needed. The algorithm is applicable to a
large class of systems, and does not rely on an explicit formula of
the large deviation rate function -- only the Hamiltonian of the
theory is needed. We derived specific reductions that are possible in
regularly occurring special cases, such as SDEs with additive or
multiplicative noise. Furthermore, we discussed optimizations for
SPDEs with additive noise and commented on how to improve numerical
efficiency. 

The performances of the new gMAM algorithm were illustrated in a
series of applications arising from different fields and involving
different types of models, like S(P)DEs with additive and
multiplicative Gaussian noises, Markov jump processes, or slow-fast
systems.

\subsection*{Acknowledgement}
We would like to thank M. Cates, A. Donev, and F. Bouchet for
helpful discussions. This research was supported in part by the NSF
grants DMR-1207432 (T. Grafke), DMS-1108780 and DMS-1522737
(T. Sch\"afer), and  DMS-1522767 (E. Vanden-Eijnden).

\bibliographystyle{plainnat}
\bibliography{arxiv}

\end{document}